\documentclass{amsart}
\usepackage{amscd, amssymb, portland, calc, ifthen}
\usepackage[matrix,arrow,curve,cmtip]{xy}

\usepackage{eucal}
\usepackage{mathrsfs}
\DeclareMathAlphabet{\mathpzc}{OT1}{pzc}{m}{it}

\newcommand{\ep}{\varepsilon}

\newcommand{\mathbold}{\bf}
\newcommand{\pser}[1]{{[\![{#1}]\!]}}

\newcommand{\longlabelmap}[1]{{\,\buildrel #1\over\longrightarrow\,}}

\newcommand{\longmap}{{\,\longrightarrow\,}}

\newcommand{\Spec}{{\mathrm{Spec}}}
\newcommand{\Mod}{{\:\mathrm{mod}\:}}

\newcommand{\sO}{{\mathcal{O}}}

\newcommand{\g}{{\mathfrak{g}}}

\newcommand{\ff}{{\mathbold F}}
\newcommand{\qq}{{\mathbold Q}}

\newcommand{\zz}{{\mathbold Z}}

\newcommand{\nn}{{\mathbold N}}

\newcommand{\m}{{\mathfrak{m}}}

\newcommand{\comment}[1]{}

\renewcommand{\geq}{\geqslant}

\let \: \colon
\newcommand{\id}{{\mathrm{id}}}

\newcommand{\Ring}{\mathsf{Ring}}

\newcommand{\Alg}[1]{\mathsf{Ring}_{{#1}}}
\renewcommand{\Mod}{\mathsf{Mod}}

\newcommand{\ct}[1]{{C_{#1}}}

\newcommand{\bcp}{\odot} 
\newcommand{\lcp}{\odot} 
\newcommand{\swe}{\circledast}
\newcommand{\vbl}{-}
\newcommand{\br}[1]{{\mathsf{BR}_{{#1}}}}
\newcommand{\unt}[1]{{{#1}\langle e\rangle}}
\newcommand{\ab}[1]{{\langle{#1}\rangle}}
\newcommand{\sd}[1]{\rtimes_{{#1}}}
\newcommand{\tn}{\otimes}
\newcommand{\dfn}{:=}

\newtheoremstyle{mythm}{}{}%
  {\itshape}
  {}
  {\bfseries}
  {}
  { }
  {\thmnumber{#2.\hspace{1.5mm}}\thmname{#1}\thmnote{ #3}.}

\newtheoremstyle{myrmk}{}{}%
  {}
  {}
  {}
  {}
  { }
  {{\bfseries\thmnumber{#2.\hspace{1.5mm}}}{\itshape\thmname{#1}}\thmnote{ #3}.}

\numberwithin{equation}{subsection}

\theoremstyle{mythm}
\newtheorem{thm}[subsection]{Theorem}
\newtheorem{prop}[subsection]{Proposition}
\newtheorem{lemma}[subsection]{Lemma}
\newtheorem{cor}[subsection]{Corollary}

\theoremstyle{myrmk}
\newtheorem{rmk}[subsection]{Remark}

\theoremstyle{plain}
\newtheorem*{thm*}{Theorem}
\newtheorem*{prop*}{Proposition}
\newtheorem*{cor*}{Corollary}
\newtheorem*{conj*}{Conjecture}

\theoremstyle{definition}

\makeatletter
\def\@seccntformat#1{\@ifundefined{#1@cntformat}%
{\csname the#1\endcsname\quad}
{\csname #1@cntformat\endcsname}
}
\def\section@cntformat{\thesection.\enspace}
\def\subsection@cntformat{\thesubsection.}
\makeatother

\newcommand\mnote[1]{}
\newcommand\bel[1]{{\mnote{#1}}\begin{equation}\label{#1}}
\newcommand\lb[1]{\label{#1}\mnote{#1}}

\newcounter{hour}\newcounter{minute}
\newcommand{\printtime}{\setcounter{hour}{\time/60}%
	\setcounter{minute}{\time-\value{hour}*60}%
	\ifthenelse{\value{hour}<10}{0}{}\thehour:%
	\ifthenelse{\value{minute}<10}{0}{}\theminute}

\begin{document}

\title{Plethystic algebra}
\author[J.~Borger, B.~Wieland]{James Borger, Ben Wieland}
\address{University of Chicago, Chicago, Illinois, USA}
\email{borger@math.uchicago.edu, wieland@math.uchciago.edu}
\thanks{Both authors were partially supported by the National Science
Foundation}

\begin{abstract}

The notion of a $\zz$-algebra has a non-linear analogue, whose purpose
it is to control operations on commutative rings rather than linear
operations on abelian groups.  These {\em plethories} can also be
considered non-linear generalizations of cocommutative bialgebras.
We establish a number of category-theoretic facts about
plethories and their actions, including a Tannaka--Krein-style
reconstruction theorem.  We show that the classical ring of Witt
vectors, with all its concomitant structure, can be understood in a
formula-free way in terms of a plethystic version of an affine blow-up
applied to the plethory generated by the Frobenius map.  We also
discuss the linear and infinitesimal structure of plethories and
explain how this gives Bloch's Frobenius operator on the
de~Rham--Witt complex.
\end{abstract}

\maketitle
\section*{Introduction}

Consider an example from arithmetic.
Let $p$ be a prime number.  Recall that for (commutative) rings $R$, the
ring $W(R)$ of ($p$-typical) Witt vectors is usually defined to be
the unique ring structure on the set $R^\nn$ which is functorial in
$R$ and such that the map
\[
(r_0,r_1,\dots) \mapsto (r_0,r_0^p+pr_1, r_0^{p^2}+pr_1^p+p^2r_2,\dots)
\]
is a ring homomorphism, the target having the usual, product ring
structure.  If $R$ is a perfect field of characteristic $p$, then
$W(R)$ is the unique complete discrete valuation ring whose maximal
ideal is generated by $p$ and whose residue field is $R$.  However, in
almost all other cases, $W(R)$ is pathological by the usual
standards of commutative algebra.  For example, $W(\ff_p[x])$ is not
noetherian.

It is nevertheless an established fact that $W(R)$ is an important
object.  For example, if $R$ is the coordinate ring of a smooth affine
variety over a perfect field of characteristic $p$, there is a certain
quotient of the de~Rham complex of $W(R)$, called the de~Rham--Witt
complex of $R$, whose cohomology is naturally the crystalline
cohomology of $R$.  But it is not at all clear from the definition
above what the proper way to think about $W(R)$ is, much less why it
is even reasonable to consider it in the first place.  The presence of
certain natural structure, for example, a multiplicative map $R\to
W(R)$ and a ring map $W(R)\to W(W(R))$ adds to the mystery.  And so we
have a question: is there a definition given purely in terms of
algebraic structure rather than somewhat mysterious formulas, and
is there a point of view from which this definition will be seen as
routine and not the result of some intangible inspiration?

The purpose of this paper is to discuss an algebraic theory of which a
particular instance gives a formal answer to these questions and to
write down some basic definitions and facts.  For any
(commutative) ring $k$, we define a $k$-plethory to be a commutative
$k$-algebra together with a comonad structure on the covariant functor
it represents, much as a $k$-algebra is the same as a $k$-module that
represents a comonad.  So, just as a $k$-algebra is exactly the
structure that knows how to act on a $k$-module, a $k$-plethory is the
structure that knows how to act on a commutative $k$-algebra.  It is
not so surprising that this analogy extends further:

\begin{table}[h]
\begin{tabular}{|c|c|}
\hline
Linear $/k$ & Non-linear $/k$ \\
\hline \hline 
$k$-modules $M$         & Commutative $k$-algebras $R$ \\
$k$-$k$-bimodules $N$   & $k$-$k$-birings $S$ \\
$\text{Hom}_k(N,M)$	& $\text{Hom}_{k{\text{-alg}}}(S,R)$ \\
$N\tn_k M$		& $S\bcp_k R$ \\
$k=\tn$-unit 		& $k[e]=\bcp$-unit \\
$k$-algebras $A$ 	& $k$-plethories $P$ \\
$A$-modules		& $P$-rings \\
$A$-$A'$-bimodules	& $P$-$P'$-birings \\
\hline
\end{tabular}
\end{table}

This is explained in section one.  In fact, as George Bergman has
informed us, this picture has been known in the
universal-algebra community, under quite similar terminology and
notation, since Tall and Wraith's paper~\cite{Tall-Wraith} in 1970.
(See also~\cite{Wraith},\cite{Bergman}.)
For those familiar with their work,
parts of the first sections will be very familiar.

The description of the ring of Witt vectors from this point of view
is that there is a $\zz$-plethory $\Lambda_p$, 
and $W(R)$ is simply the $\Lambda_p$-ring co-induced from
the ring $R$ (which observation allows us to define a Witt
ring for any plethory), and so the only thing left is to give a natural
construction of $\Lambda_p$.  This is done by a process we call {\em
amplification} and which is formally similar to performing an affine
blow-up in commutative algebra.  We will give some idea of this
procedure below.


In section two, we give some examples of plethories.  The most basic
is the symmetric algebra $S(A)$ of any cocommutative bialgebra $A$; in
particular, if $A$ is a group algebra $\zz G$, then $S(A)$ is
the free polynomial algebra on the set underlying $G$.  These
plethories are less interesting because their actions on rings
can be described entirely in terms of the original bialgebra $A$; for
example, an action of the plethory $S(\zz G)$ 
is the same as an action
of the group $G$.  But even in this case, there can be more maps between
two such plethories than there are between the bialgebras, and in some
sense, this is ultimately responsible for existence of $\Lambda_p$ and
hence the $p$-typical Witt ring.

The ring $\Lambda$ of symmetric functions in infinitely many variables
is a better example.  The composition law of $\Lambda$ is given by the
operation known as plethysm in the theory of symmetric functions and
is what gives plethories their name.  An action of $\Lambda$ on a ring
$R$ is the same as a $\lambda$-ring structure on $R$, and in
contrast to plethories of the form $S(A)$, a $\Lambda$-action cannot
in general be described in terms of a bialgebra action.  We also give
an explicit description of $\Lambda_p$, the plethory responsible for
the $p$-typical Witt ring, in terms of symmetric functions.  Of
course, this description is really quite close to a standard treatment
of the Witt ring and is still a bit unsatisfying.  In section three,
we give explicit examples of $P$-Witt rings for various plethories
$P$.

In section four, we discuss the restriction, induction, and
co-induction functors for a morphism $P\to Q$ of plethories, and we
state the reconstruction theorem.  As always, the content of such a
theorem is entirely category theoretic (Beck's theorem). 
All the same, the result is worth stating:

\begin{thm*}
Let $\mathsf{C}$ be a category that has all limits and colimits, let
$U$ be a functor from $\mathsf{C}$ to the category of rings.  If $U$
has both a left and a right adjoint and has the property that a map
$f$ in $\mathsf{C}$ is an isomorphism if $U(f)$ is, then $\mathsf{C}$
is the category of $P$-rings for a unique $k$-plethory $P$, and under
this identification, $U$ is the forgetful functor from $P$-rings to
rings.
\end{thm*}


In section seven, we explain amplification, the blow-up-like process
we mentioned above.  Let $\sO$ be a Dedekind
domain, for example the ring of integers in a local or global field
or the coordinate ring of a smooth curve.
Let $\m$ be an ideal in $\sO$, let $P$ be an $\sO$-plethory, let $Q$
be an $\sO/\m$-plethory, and let $P\to Q$ be a surjective map of
plethories.  We say a $P$-ring $R$ is a {\em $P$-deformation of a
$Q$-ring} if it is $\m$-torsion-free and the action of $P$ on $R/\m R$
factors through the map $P\to Q$.  

\begin{thm*}
There is an $\sO$-plethory $P'$ that is universal among those that are
equipped with a map from $P$ making them $P$-deformations of
$Q$-rings.  Furthermore, $P'$ has the property that $P$-deformations of
$Q$-rings are the same as $P'$-rings that are $\m$-torsion-free.
\end{thm*}

We say $P'$ is the amplification of $P$ along $Q$.

In sections eight through eleven, we define what could be called the
linearization of a plethory $P$.  It involves two
structures: $A_P$, the set of elements of $P$ that act
additively on any $P$-ring, and $\ct{P}$, the cotangent space to
the spectrum of $P$ at $0$.

\begin{thm*}
Both $A_P$ and $\ct{P}$ are (generally non-commutative) 
algebras equipped with maps from $k$, and
under certain flatness or splitting hypotheses, the following hold:
$A_P$ is a cocommutative twisted $k$-bialgebra, there is a coaction of
$A_P$ on the algebra $\ct{P}$, and the map $A_P\to\ct{P}$ is
$A$-coequivariant.
\end{thm*}

We stop short of investigating representations of such linear
structures. 

If $R\to R'$ is a map of $P$-rings with kernel $I$, then all that
remains on the conormal module $I/I^2$ of the action of $P$ is an
action of $\ct{P}$.  In particular, $\ct{P}$ acts on the K\"ahler
differentials of any $P$-ring.  In the special case when $P=\Lambda_p$
and $R=W(S)$, for some ring $S$, 
this additional structure is essentially a lift of Bloch's Frobenius
operator on the de~Rham--Witt complex.

The final section of the paper is the reason why the others exist, and
we encourage the reader to look at it first.  Here, we consider
$\Lambda_p$ and other classical constructions from the point of view
of the general theory.  For example, we give a satisfying construction
of $\Lambda_p$: Let $\unt{\ff_p}$ be the trivial
$\ff_p$-plethory; its bialgebra of additive elements has a canonical
deformation to a $\zz$-bialgebra, and let $P$ be the free
$\zz$-plethory on this.  Then $\Lambda_p$ is the amplification of $P$
along $\unt{\ff_p}$.  
Essentially the same procedure,
applied to rings of integers in general number fields, gives at once ramified
and twisted generalizations. 

An action of this amplification on a $p$-torsion-free ring $R$ is,
essentially by definition, the same as a lift of the Frobenius
endomorphism of $R/pR$.  The content of the statement that
the $\Lambda_p$-ring co-induced by $R$
agrees with the classical $W(R)$ is ultimately just
Cartier's Dieudonn\'e--Dwork lemma.  Thus it would be accurate to view
amplifications as the framework where Joyal's approach to the classical
Witt vectors~\cite{Joyal:Witt} naturally lives.

The last section also has explicit descriptions of the
linearizations of $\Lambda_p$, $\Lambda$, and similar plethories.


On a final note, this paper does not even contain the basics of the
theory, and there are still many simple mysteries.  For example, the
existence of non-linear plethories, those that do not come from
(possibly twisted) bialgebras, may be a purely arithmetic phenomenon:
we know of no non-linear plethory over a $\qq$-algebra.  For a broader
example, the category of $P$-rings is, on the one hand, a
generalization of the category of rings and, on the other, an analogue
of the category of modules over an algebra.  And so it is natural to
ask which notions in commutative algebra and algebraic geometry can be
generalized to $P$-rings for general $P$ and, in the other direction,
which notions in the theory of modules over algebras have analogues in
the theory actions of plethories on rings.  It would be quite
interesting to see how far these analogies can be taken.


It is a pleasure to thank Spencer Bloch and David Ben-Zvi for their
thoughts, both mathematical and terminological.

\setcounter{section}{0}
\section*{Conventions}
\lb{sec-notation}

The word {\em ring} is short for commutative ring, but we make no
commutativity restriction on the word {\em algebra}.  A $k$-ring is
then a commutative $k$-algebra.  All these objects are assumed to be
associative and unital, and all morphisms are unital.

We use the language of coalgebras
extensively; D{\u{a}}sc{\u{a}}lescu, N{\u{a}}st{\u{a}}sescu,
and Raianu's book~\cite{Dascalescu:HA} is more than enough.

For categorical terminology, we refer to Mac Lane's book~\cite{MacLane:CWM}.
In particular, we find it convenient to write
$\mathsf{C}(X,Y)$ for the set of morphisms between objects $X$ and $Y$
of a category $\mathsf{C}$.

$\nn$ denotes the set $\{0,1,2,\dots\}$.

\section{Plethories and the composition product}

Let $k,k',k''$ be rings. 

A $k$-$k$-biring is a $k$-ring that represents a
functor $\Ring_k\to\Ring_k$.  Composition of such functors yields a
monoidal structure on the category of $k$-$k$-birings.  We then define
a $k$-plethory to be a monoid in this category, much as one could
define a $k$-algebra to be a monoid in the category of
$k$-$k$-bimodules.  Finally, the category of $k$-$k$-birings acts on
the category of $k$-rings, and we define a
$P$-ring to be a ring together with an action of the $k$-plethory $P$.

We spell this out in some detail and give a number of immediate
consequences of the definitions.  We also give many examples in this
section, but they are all trivial, and so 
the reader may want to look ahead at the more interesting
examples in sections two and three.



\subsection{}
\lb{sbsc-biring}
A $k$-$k'$-{\em biring} is a $k$-ring $S$, together with a lift of the
covariant functor it represents to a functor
$\Ring_k\to\Ring_{k'}$. Equivalently, it is the structure on $S$ of a
$k'$-ring object in the opposite category of $\Ring_k$. Or in
Grothendieck's terminology, this is the structure on $\Spec\,S$ of a
commutative $k'$-algebra scheme over $\Spec\,k$. Explicitly, $S$ is a
$k$-ring with the following additional maps (all of $k$-rings except
(3)):
 
\begin{enumerate}
\item {\em coaddition:} \lb{axiom-coaddition}
a cocommutative coassociative map $\Delta^+\:S\to
S\tn_{k} S$ for which there exists
a counit $\ep^+\:S\to k$ and an antipode $\sigma\:S\to S$,
\item {\em comultiplication:} \lb{axiom-comultiplication}
a cocommutative coassociative map 
$\Delta^\times\: S\to S\tn_{k} S$ which codistributes over
$\Delta^+$ and for which there exists a counit $\ep^\times\:S\to k$,
\item {\em co-$k'$-linear structure:} \lb{axiom-unit}
a map $\beta\:k'\to\Alg{k}(S,k)$
of rings, where the ring structure on $\Alg{k}(S,k)$ is given
by~(\ref{axiom-coaddition})
and~(\ref{axiom-comultiplication}).
\end{enumerate}
Note that, as usual, $\ep^+$, $\sigma$, and $\ep^\times$ are unique if
they exist.  Also note that omitting axiom~(\ref{axiom-unit}) leaves us
with the notion of $k$-$\zz$-biring.  Finally, in the case of
$k$-plethories, we will take $k=k'$, but at this point it is best to
keep the roles separate.

A morphism of $k$-$k'$-birings is a map of $k$-rings which
preserves all the structure above.  The category of $k$-$k'$-birings
is denoted $\br{k,k'}$.  Given a map $k''\to k'$, we can view a
$S$ as a $k$-$k''$-biring, which we still denote $S$, somewhat abusively.

Let $\ell$ and $\ell'$ be rings, and let $T$ be a
$\ell$-$\ell'$-biring.  A morphism $S\to T$ of birings is
the following data: a ring map $k\to\ell$, a ring map $k'\to
\ell'$, and a map $\ell\tn_k S \to T$ of
$\ell$-$k'$-birings.  The category of birings is denoted $\br{}$.
When necessary, we will distinguish the structure maps of birings by using
subscripts: $\Delta^+_S, \ep^\times_S,$ and so on.  We will also often
use without comment the notation $\Delta^+p=\sum_i p^{(1)}_i\tn
p^{(2)}_i$ and $\Delta^\times p=\sum_i p^{[1]}_i\tn p^{[2]}_i$.


\subsection{} \lb{sbsc-ident}
{\em Examples.}
\begin{enumerate}
\item $k$ itself is the initial $k$-$k'$-biring, representing the 
constant functor giving the zero ring.
\item Let $\unt{k}$ denote the $k$-$k$-biring that represents the identity
functor on $\Alg{k}$.  Thus $\unt{k}$ is canonically the ring $k[e]$ with
$\Delta^+(e)=e\tn 1 + 1\tn e,\Delta^\times(e)=e\tn e,
\beta(c)(e)=c$ (and $\ep^+(e)=0, \ep^\times(e)=1, \sigma(e)=-e$).
\item
If $k'$ is finite, then the collection of set maps $k'\to k$ is
naturally a $k$-$k'$-biring.  The
$k$-ring structure is given by pointwise addition and
multiplication, and the coring structure is given by the ring
structure on $k'$.  For example, $\Delta^+$ is the composite
$k^{k'}\to k^{k'\times k'} = k^{k'}\otimes k^{k'}$, where the first
map is given by addition on $k'$. If $k'$ is not finite, there are
topological issues, which could surely be avoided by considering
pro-representable functors from $\Alg{k}$ to $\Alg{k'}$.
\end{enumerate}

\vspace{2ex}

Recall that the action of a $k$-algebra $A$ on a $k$-module $M$ can be given in
three ways: as a map $A\otimes_k M\to M$, as a map $M\to \Mod_k(A,M)$,
or as a map $A\to \Mod_k(M,M)$. In fact, we have the same choices when
defining the multiplication map on $A$ itself.
The Witt vector approach to operations on rings follows the second,
comonadic model, but we will follow the first, monadic one.
The third approach encounters the topological problems mentioned in
the example above.

We now define the analogue of the tensor product.


\subsection{} \lb{sbsc-comp-prod}
{\em Functor $\vbl\bcp_{k'}\vbl\:\br{k,k'}\times\Alg{k'}\to\Alg{k}$.} 
Take
$S\in\br{k,k'}$
and
$R\in\Ring_{k'}$.
Then $S\bcp_{k'} R$ is defined to be the $k$-ring
generated by symbols $s\lcp r$, for all $s\in S, r\in R$, subject to
the relations (for all $s,s'\in S, r,r'\in R, c\in k'$)
\bel{equ-alg-rel}
ss'\lcp r = (s\lcp r)(s'\lcp r),\quad
(s+s')\lcp r =(s\lcp r) + (s'\lcp r),\quad c\lcp r = c
\end{equation}
and 
\bel{equ-coalg-rel}
\begin{split}
s\lcp(r+r')&=\Delta_S^+(s)(r,r')\dfn\sum_i (s_i^{(1)}\lcp r)(s_i^{(2)}\lcp r'), \\
s\lcp(rr')&=\Delta_S^\times(s)(r,r')\dfn\sum_i (s_i^{[1]}\lcp r)(s_i^{[2]}\lcp r'), \\
s\lcp c &= \beta(c)(s).
\end{split}
\end{equation}
This operation is called the  {\em composition product}
and is clearly functorial in both $R$ and $S$.

As in linear algebra, where a tensor
$a\tn b$ reminds us of the formal composition of operators $a$ and
$b$ or the formal evaluation of an operator $a$ at $b$, the symbol $s\lcp
r$ is intended to remind us of the composition $s\circ r$ of possibly
non-linear functions or the formal evaluation of a function $s$ at
$r$.  Thus the meaning of (\ref{equ-alg-rel}) is that ring operations
on functions are defined pointwise, and the meaning
of (\ref{equ-coalg-rel}) is that there is extra structure on our ring of
functions that controls how they respect sums, products, and constant
functions. For example, if $S$ is the biring of~\ref{sbsc-ident}(3),
the evaluation map $S\bcp_{k'}k' \to k$ given by $s\lcp r\mapsto s(r)$
is a well-defined ring map.

\begin{prop}
\lb{prop-comp-prod-adj}
Let $S$ be a $k$-$k'$-biring.  The functor $S\bcp_{k'}\vbl$ 
is the left adjoint of $\Alg{k}(S,\vbl)$. 
\end{prop}

In other words, for $R_1\in\Alg{k}$, $R_2\in\Alg{k'}$
we have
\[
\Alg{k}(S\bcp_{k'}R_2,R_1)=\Alg{k'}(R_2,\Alg{k}(S,R_1)).
\]
The proof is completely straightforward.
We leave it, as well as the task of specifying the unit and counit
of the adjunction, to the reader.

\subsection{} \lb{examples-comp-prod}
{\em Examples.} 
\begin{enumerate}
\item
There are natural identifications $S\bcp_{k'} \unt{k'}=S$, 
$\unt{k'} \bcp_{k'} R=R$,
$S\bcp_{k'} k' = k$, and $k\bcp_{k'} R = k$.
\item
If $k'\to \ell'$ is a ring map, then $\unt{\ell'}\bcp_{k'}
R=\ell'\tn_{k'} R$.
\item 
$k$-$\ell'$-biring structures on $S$ compatible with the given
$k$-$k'$-biring structure are the same, under adjunction, as maps
$S\bcp_{k'}\ell'\to k$ of $k$-rings.  
\item
If $k\to \ell$ is a ring map, we
have $(\ell\tn_k S)\bcp_{k'} R = \ell\tn_k (S\bcp_{k'} R)$.
\item
The composition product distributes over arbitrary tensor products:
\begin{align*}
\Big(\bigotimes S_i\Big)\bcp_{k'} R &= \bigotimes \big(S_i\bcp_{k'}R\big), \\
S\bcp_{k'}\Big(\bigotimes R_i\Big) &= \bigotimes \big(S\bcp_{k'}R_i\big).
\end{align*}
\end{enumerate}

\subsection{}
\lb{sbsc-biring-composition}
If $R$ is not only a $k'$-ring but a $k'$-$k''$-biring, then the
functor 
\[
\Alg{k}(S\bcp_{k'}R,\vbl)=\Alg{k'}(R,\Alg{k}(S,\vbl))
\]
naturally takes values in $k''$-rings, and so $S\bcp_{k'}R$ is
naturally a $k$-$k''$-biring.  One can also see this directly in terms of the
structure maps $\Delta^+$ and so on by using the fact that the composition
product distributes over tensor products.
If $k=k'=k''$, the composition product gives a
monoidal structure on the category of $k$-$k$-birings with unit
$\unt{k}=k[e]$ of~\ref{sbsc-ident}.  As is generally true
with composition or the
tensor product of bimodules, this monoidal structure not
symmetric.


\begin{rmk}
Note that, in contrast to the analogous statement for bimodules, it is
generally not true that a $k$-$k''$-biring structure on $R$ induces
$k'$-$k''$-biring structure on the $k$-ring $\Alg{k}(S,R)$.
\end{rmk}

\subsection{}
\lb{sbsc-plethory-def}
A {\em $k$-plethory} is a monoid in the category of $k$-$k$-birings,
that is, it is a biring $P$ equipped with an associative map of
birings $\circ\:P\bcp_k P\to P$ and unit $\unt{k}\to P$.  For example,
$\unt{k}=k[e]$ with $\circ$ taken as in~\ref{examples-comp-prod}(1)
(that is, composition of polynomials) is a
$k$-plethory.  The image of
$e$ under the unit map $\unt{k}\to P$ is denoted $e$ (or $e_P$); together
with $\circ$, it gives the set underlying $P$ a monoid structure. The
ring $k$ is called the {\em ring of scalars} of $P$.  

If $P'$ is a $k'$-plethory,
a morphism $P\to P'$ of plethories is a morphism $k\to k'$ plus
a morphism $\varphi\:P\to P'$ of birings which is also a morphism of monoids.
This is equivalent to requiring that
\[
\xymatrix{
\unt{k'}\bcp_k P\bcp_k P \ar^-{\varphi\lcp 1}[r]\ar^-{1\lcp\circ}[dd]
  & P'\bcp_k P\ar@{=}[r]
  & P'\bcp_{k'}\unt{k'}\bcp_k P \ar^-{1\lcp\varphi}[d] \\
& & P'\bcp_{k'} P' \ar^-{\circ}[d] \\
\unt{k'}\bcp_k P \ar@{=}[r] & k'\tn_k P\ar^{\varphi}[r] & P'
}
\]
be a commutative diagram of $k'$-$k$-birings.  If $k=k'$, the diagram
simplifies to the obvious one.
If we are already given a map $k\to k'$, then we will always assume
the map of scalars is the same as the given map.
It is easy to see that $\unt{k}$ is the initial
$k$-plethory and $\unt{\zz}$ is the initial plethory.

\subsection{} \lb{sbsc-P-action}
A {\em (left) action} of $P$ on a $k$-ring $R$ is defined as usual in
the theory of monoidal categories; in this case it means a map
$\circ\:P\bcp R\to R$ such that $(\alpha\circ\beta)\circ r =
\alpha\circ(\beta\circ r)$ and $e\circ r = r$ for all $\alpha,\beta\in
P, r\in R$.  We also denote $\alpha\circ r$ by $\alpha(r)$.  A {\em
$P$-ring} is a $k$-ring equipped with an action of $P$. (There is no
danger of a conflict in terminology with a ring equipped with a ring
map from $P$ because we never use such structures in this paper.)  A
morphism of $P$-rings is a map of rings that makes the obvious diagram
commute; equivalently, it is a map of rings that is $P$-equivariant as
a map of sets acted on by the monoid ($P,\circ)$.
The category of $P$-rings is denoted $\Ring_P$.  

If $S$ is a $k$-$k'$-biring, we say $P$ {\em acts on $S$ as a
$k$-$k'$-biring} if $\circ\:P\bcp S\to S$ is a map of
$k$-$k'$-birings.  Such an action is the same as a functorial
collection of $k'$-ring structures on the sets $\Ring_P(S,R)$ such
that the maps $\Ring_P(S,R)\hookrightarrow\Alg{k}(S,R)$ are maps of
$k'$-rings.

A {\em right action} of a $k'$-plethory $P'$ on a $k$-$k'$-biring is a
map $\circ\:R\bcp_{k'} P'\to R$ of $k$-$k'$-birings compatible with
$\circ$ and $e$ in the obvious way.  A map of right $P'$-rings is
$P'$-equivariant map of $k$-$k'$-birings.  A $P$-$P'$-biring is a
$k$-$k'$-biring equipped with a left action of $P$ as a $k$-$k'$
biring and a commuting right action of $P'$. The category of
$P$-$P'$-birings is denoted $\br{P,P'}$, morphisms being maps of
birings that are both $P$-equivariant and $P'$-equivariant.
A $P$-$P'$-biring is the same as a represented functor $\Ring_P\to\Ring_{P'}$.

\subsection{} \lb{sbsc-F-W}
A $k$-plethory structure on a $k$-$k$-biring $P$ is the same as a
monad structure on the functor $P\bcp_k\vbl$ and, by adjunction, also
the same as a comonad structure on the functor $\Alg{k}(P,\vbl)$. An
action of $P$ on $R$ is the same as the structure on $R$ of an algebra
over the monad or a coalgebra over the comonad.

Thus $\Ring_{P}$ has all limits and colimits, the forgetful functor
$U\:\Ring_P\to\Alg{k}$ preserves them, and the functors
$P\bcp_{k}\vbl$ and $\Alg{k}(P,\vbl)$ lift to give left and,
respectively, right adjoints to $U$.  (These functors could well be called
restriction, induction, and co-induction for the map $\unt{k}\to P$.
We postpone the treatment of these functors for general maps of
plethories until section four.)  In particular, the underlying
$k$-ring of a (co)limit of $P$-rings is the (co)limit in that
category and there exists a unique compatible $P$-ring structure on
it.  We give a converse to all this in
section~\ref{sec-reconstruction}.

We often denote the functor $\Alg{k}(P,\vbl)$ by $W_P(\vbl)$ and
call the $P$-ring $W_P(R)$ the {\em $P$-Witt ring} of $R$.  The reason for
this terminology will be made clear in section~\ref{sec-Witt}.

\subsection{} \lb{examples-P-rings}
{\em Examples.}
\begin{enumerate}
\item
If $k$ is finite, the biring of set maps $k\to k$ is a
$k$-plethory, with $\circ$ given by composition of functions.  In
particular, $0$ is a plethory over the ring $0$.  It is the terminal
plethory, and of course the only $0$-ring is $0$.
\item \lb{example-P-is-free}
A plethory $P$ clearly acts on itself on the left (and also the
right).  It is in fact the free $P$-ring on one element: morphisms in
$\Ring_P$ from $P$ to another object are the same as elements of the
underlying ring, a map $\varphi\:P\to R$ corresponding to the element
$\varphi(e)$ 
in $R$, and an element $r\in R$ corresponding to the map
$\alpha\mapsto\alpha(r)$. 
The morphisms $P\to k$ corresponding to $r=0$ and $r=1$ are
$\ep^+$ and $\ep^\times$.  More generally, the morphism $P\to k$
corresponding to $c\in k$ is $\beta(c)$.
\item
\lb{example-counits}
The identification $P\bcp_k k = k$ is an action of $P$ on $k$, and if
$R$ is any $P$-ring, the structure map $k\to R$ is a map of $P$-rings
simply by the third relation of (\ref{equ-coalg-rel}).  Therefore, $k$ is
the initial $P$-ring.  Similarly, the identification $k\bcp_k P=k$
gives $k$ the structure of a $P$-$P$-biring, and it is the initial
$P$-$P$-biring.
\item
If $k'$ is a $P$-ring, the natural $k'$-map 
\[
(k'\tn_k P)\bcp_k k' = k'\tn_k (P\bcp_k k') \to k' 
\]
gives (by~\ref{examples-comp-prod})
$k'\tn_k P$ the structure of a $k'$-$k'$-biring.  We will
see below that $k'\tn_k P$ even has a natural $k'$-plethory structure.
\end{enumerate}

\begin{prop}
\lb{prop-P-equ}
Let $P$ be a $k$-plethory.  Then
the $k$-ring morphisms $\Delta^+_P$, $\Delta^\times_P$, $\ep^+_P$, and
$\ep^\times_P$ are in fact $P$-ring morphisms. For any $A\in \Ring_P$,
the unit $\eta_A\:k\to A$ and multiplication $m_A\:A\tn_k A\to A$
are $P$-ring morphisms.
\end{prop}
\begin{proof}
The unit and counits were discussed in
\ref{examples-P-rings} (\ref{example-counits}) and (\ref{example-P-is-free}).
Multiplication is the
coproduct of the identity with itself.

By~\ref{examples-P-rings}(\ref{example-P-is-free}), the $P$-ring $P$
represents the forgetful functor $U'$ from $\Ring_P$ to the category
of sets and $P\tn_k P$ represents the functor $U'\times U'$.  But
these factor through the category of rings, and so there are
natural transformations $U'\times U'\to U'$, one for addition and one
for multiplication. Thus there are maps $P\to P\tn_k P$ in
$\Ring_P$. The one for addition is the map that sends $e$ to $1\tn
e+e\tn1$, and thus sends $\alpha$ to $\Delta^+(\alpha)(1\tn e,e\tn
1)=\Delta^+(\alpha)$.  Similarly, the one for multiplication is
$\Delta^\times$.
\end{proof}

\subsection{} \lb{sbsc-base-change}
{\em Base change of plethories.}
If $k'$ is a $P$-ring, then the $k'$-$k$-biring $k'\tn_k P$ has a
$k'$-$k'$-biring structure (\ref{examples-P-rings}).
Even further, the $k'$-ring map (using~\ref{examples-comp-prod}(4))
\[
(k'\tn_k P)\bcp_k (k'\tn_k P) = k'\tn_k (P\bcp_k
(k'\tn_k P)) \longlabelmap{1\tn\circ} k'\tn_k (k'\tn_k P)\longmap k'\tn_k P
\]
descends to a map
\[
(k'\tn_k P)\bcp_{k'} (k'\tn_k P) \longmap k'\tn_k P,
\]
which gives $k'\tn_k P$ the structure of a $k'$-plethory.

Conversely,
if $k'\tn P$ is a $k'$-plethory, then $P$ acts on $k'$ by way of 
$k'\tn P$.  Note that not only does the plethory
structure on $k'\tn P$ depend on the action of $P$ on $k'$, 
there may not exist even one such action.
For example, there is no action of the $\zz$-plethory $\Lambda_p$
(of~\ref{sbsc-lambda-p}) on $\ff_p$.

We leave it as an exercise to show that a $k'\tn P$-action on a
$k'$-ring $R$ is the same as a $P$-action on the underlying $k$-ring
compatible with the given action on $k'$.

\section{Examples of plethories}

Before continuing with the theory, let us give some basic examples of
plethories.

\subsection{} 
\lb{sbsc-free-on-biring}
{\em Free plethory on a biring.}
Let $k$ be a ring, and let $S$ be a $k$-$k$-biring.  There is a
plethystic analogue of the tensor algebra: a $k$-plethory $Q$, with a
$k$-$k$-biring map $S\to Q$, which is initial in the category of
such plethories.

Put
$$
Q = \bigotimes_{n\ge0} S^{\lcp n}.
$$
The system of maps
\begin{align*}
S^{\lcp i} \bcp S^{\lcp j} &\longmap S^{\lcp (i+j)} \\
(s_1\lcp\cdots\lcp s_i) \lcp (t_1\lcp\cdots\lcp t_j)&\mapsto
s_1\lcp\cdots\lcp s_i \lcp t_1\lcp\cdots\lcp t_j
\end{align*}
induces a map
\[
Q\bcp Q = \bigotimes_{i,j} S^{\lcp i} \bcp S^{\lcp j} \longmap 
\bigotimes_n S^{\lcp n} = Q,
\]
which is clearly associative.  This gives $Q$ the structure of a
$k$-plethory with a map $\unt{k}=S^{\bcp 0}\to Q$ of $k$-plethories.

A $Q$-action on a ring $R$ is then the same as a map $S\bcp R\to R$ of
rings.

\subsection{} \lb{sbsc-free-on-bialg}
{\em Free plethory on a cocommutative bialgebra.}
First, let $A$ be a cocommutative coalgebra over $k$; denote its
comultiplication map by $\Delta$ and its counit by $\ep$.  The
symmetric algebra $S(A)$ of $A$, viewed as a $k$-module, is of course
a $k$-ring, but the following gives it the structure of a
$k$-$k$-biring:

{\em Coadditive structure:}
The coaddition map $\Delta^+$ is the one induced by the linear map
\[
A \longmap S(A)\tn S(A),\quad a\mapsto a\tn 1 + 1\tn a
\]
The additive counit $\ep^+\:S(A)\to k$ is the
map induced by the zero map $A\to k$.

{\em Comultiplicative structure:}
$\Delta^\times$ is the map induced by the linear map
\[
A\longlabelmap{\Delta} A\tn A \longmap S(A)\tn S(A),
\]
where the right map
is the tensor square of the canonical inclusion.  The multiplicative
counit $\ep^{\times}\:S(A)\to k$ is the composite map
\[
S(A) \xrightarrow{S(\ep)} S(k) = \unt{k} \xrightarrow{\ep^\times_{\unt{k}}}
k.
\]

{\em Co-$k$-linear structure:}
The map
\[
S(A)\bcp_{\zz} k \longmap \unt{k}\bcp_{\zz} k \longmap \unt{k}\bcp_k k
= k
\]
gives $S(A)$ a $k$-$k$-biring structure by \ref{examples-comp-prod}.

\subsection{} 
{\em Isomorphism $S(A)\bcp S(B)\to S(A\tn B)$ of $k$-$k$-birings.}
Let $B$ be another cocommutative $k$-coalgebra, and let $R$ be a
$k$-ring.  Then we have
\begin{align*}
\Ring_k(S(A)\bcp S(B),R) &= \Ring_k(S(B),\Ring_k(S(A),R)) \\
			 &= \Mod_k(B,\Mod_k(A,R)) = \Mod_k(A\tn B,R)\\
			 &= \Ring_k(S(A\tn B),R)
\end{align*}
and hence a natural isomorphism $S(A)\bcp S(B)\cong S(A\tn B)$ of
$k$-rings.  Explicitly, $[a]\lcp [b]$
corresponds to $[a\tn b]$, where $[a]$ denotes the image of $a$
under the natural inclusion $A\to S(A)$ and likewise for $[b]$.
We leave the task of showing this is a map of $k$-$k$-birings to the reader.

\subsection{}
It follows that the comultiplication and the counit induce maps
\begin{gather*}
S(A) \longmap S(A)\bcp S(A) \\
S(A) \longmap \unt{k}
\end{gather*}
that give $S(A)$ the structure of a commutative comonoid in $\br{k,k}$.

\subsection{}
Now suppose $A$ is a bialgebra, that is, $A$ is equipped with maps
\begin{gather*}
A\tn A \longmap A \\
k \longmap A
\end{gather*}
of $k$-coalgebras making $A$ a monoid in the category of
$k$-coalgebras.  By the discussion above, this makes $S(A)$ a monoid
in the category of cocommutative comonoids in $\br{k,k}$.  It is in
particular a $k$-plethory.  (It could reasonably be called a
cocommutative bimonoid in $\br{k,k}$---its additional structure is
the analogue of the structure added to an algebra to make it a
cocommutative bialgebra---but because $\bcp$ is not a symmetric
operation on all of $\br{k,k}$, this terminology could be confusing.)

\begin{rmk}
Given a $k$-ring $R$, an action of the plethory $S(A)$ on $R$ is the
same as an action of the bialgebra $A$ on $R$.  We leave the precise
formulation and proof of this to the reader.  It may be worth noting
that any $k$-ring admits an $S(A)$-action in a trivial way.  This is
true by the previous remark or by using the natural map
$S(A)\to\unt{k}$ of $k$-plethories.  It is false for general
plethories.
\end{rmk}

\subsection{} \lb{eg-bialg}
{\em Examples.}
\begin{enumerate}
\item If $A$ is the group algebra $kG$ of a group (or monoid) $G$,
then $S(A)$ is the free polynomial algebra on the set underlying $G$.
For any $g\in G$, the corresponding element in $S(A)$ is ``ring-like'':
$\Delta^+(g) = g\tn 1 + 1\tn g$ and
$\Delta^{\times}(g) = g\tn g$.  An action of the plethory $S(A)$
on a ring $R$ is the same as an action of $G$ on $R$.
\item Let $\g$ be a Lie algebra over $k$, and let $A$ be its universal
enveloping algebra.  Then for all $x\in\g$, the corresponding
element $x\in S(A)$ is ``derivation-like'': $\Delta^+(x) = x\tn 1
+ 1\tn x$ and $\Delta^{\times}(x) = x\tn e + e\tn x$.  
If $\g$ is the one-dimensional Lie algebra spanned by an element
$d$, then $S(A)=k[d^{\circ\nn}]:=k[e,d,d\circ d,\dots]$, and
$S(A)$-rings are the same as $k$-rings equipped with a derivation.
\end{enumerate}

\begin{rmk}
Because of the identification $S(A)\bcp_k S(B)\to S(A\tn B)$,
there is a natural
isomorphism $S(A)\bcp_k S(B) \to S(B)\bcp_k S(A)$ of $k$-$k$-birings
given by the
canonical interchange map on the tensor product.  Explicitly, it
exchanges $[a]\lcp [b]$ and $[b]\lcp [a]$, where $a\in A, b\in B$.
There is no functorial map $S\bcp T \to T\bcp S$ for
$k$-$k$-birings $S$ and $T$ that agrees with the previous map when $S$ and
$T$ come from bialgebras.  For example, take $S=\zz[d^{\circ\nn}]$
and $T=\Lambda_p$ below.
\end{rmk}

\subsection{} 
{\em Hopf algebras.}  An antipode $s\:A\to A$ gives a map
$S(A)\to S(A)$ of $k$-$k$-birings, making $S(A)$ what could be called
a cocommutative Hopf monoid in $\br{k,k}$.

\subsection{} 
{\em Symmetric functions and $\lambda$-rings.}
Let $\Lambda$ be the ring of symmetric functions in countably many
variables, i.e., writing $\Lambda_n$ for the sub-graded-ring of
$\zz[x_1,\dots,x_n]$ ($\deg x_i=1$) of elements invariant under the
obvious action of the $n$-th symmetric group, we let $\Lambda$ be
the inverse limit of
\[
\cdots\longmap\Lambda_n\longmap \Lambda_{n-1} \longmap\cdots .
\]
in the category of graded rings.  The map above sends $x_n$ to $0$ and
sends any other $x_i$ to $x_i$.  Of course, $\Lambda$ is the free
polynomial algebra on the elementary symmetric
functions~\cite[I.2]{MacDonald:SF}, but there are many other free
generating sets, and making this or any other particular choice would
leave us with the usual formulaic mess in the theory of
$\lambda$-rings and Witt vectors.  

The ring $\Lambda$ naturally has the structure of a plethory over
$\zz$.  Because all the structure maps are already described at
various points in the second edition of
MacDonald \cite{MacDonald:SF}, we give only the
briefest descriptions here:

{\em Coadditive structure:} \cite[I.5 ex.\ 25]{MacDonald:SF} For 
$f\in\Lambda$, consider the function
\[
\Delta^+(f) = f(x_1\tn 1, 1\tn x_1, x_2\tn 1, 1\tn x_2,\dots)
\]
in the variables $x_i\tn x_j, (i,j\geq 1)$.  It is symmetric in
both factors, and so $\Delta^+$ is a ring map $\Lambda\to
\Lambda\tn_{\zz} \Lambda$.  The counit $\ep^+\:\Lambda\to k$ sends
$f$ to $f(0,0,\dots)$. 

{\em Comultiplicative structure:} \cite[I.7 ex.\ 20]{MacDonald:SF}
Similarly, consider the function
\[
\Delta^{\times}(f) = f(\dots,x_i\tn x_j,\dots)
\]
in the variables $x_i\tn x_j$.  As before, it is symmetric in both
factors, and so $\Delta^{\times}$ is a map $\Lambda\to
\Lambda\tn_{\zz} \Lambda$.  The counit $\ep^{\times}:\Lambda\to k$ sends
$f$ to $f(1,0,0,\dots)$.

{\em Monoid structure:} \cite[I.8]{MacDonald:SF} 
For $f,g\in \Lambda$, the operation known as plethysm defines $f\circ
g$: Suppose $g$ has only non-negative coefficients, and write $g$ as a
sum of monomials with coefficient $1$ in the variables $x_i$.  Then $f\circ g$
is the symmetric function obtained by substituting these monomials
into the arguments $x_1, x_2,\cdots$ of $f$.  This gives a monoid
structure with identity $x_1+x_2+\cdots$ on the set of elements with
non-negative coefficients, and this extends to a unique $\zz$-plethory
structure on all of $\Lambda$.

\begin{rmk}
\lb{rmk-lambda-basis}
By the theorem of elementary symmetric 
functions~\cite[I 2.4]{MacDonald:SF}, we have
\[
\Lambda = \zz[\lambda_1,\lambda_2,\dots],
\]
where $\lambda_1=x_1+x_2+\cdots,\lambda_2 =
x_1x_2+x_1x_3+x_2x_3+\cdots,\dots$ are the elementary symmetric
functions.  Any $\Lambda$-ring $R$ therefore has unary operations
$\lambda_1,\lambda_2,\dots$.  It is an exercise in definitions to show
that in this way, a $\Lambda$-ring structure on a ring $R$ is the same
as a $\lambda$-ring structure (which, in Grothendieck's original
terminology~\cite{AT:Lambda}, is called a special $\lambda$-ring
structure).  This was in fact one of the principal examples in Tall
and Wraith's paper~\cite{Tall-Wraith}.

Let $\psi_n$ denote the $n$-th Adams operation:
\[
\psi_n = x_1^n + x_2^n + \cdots.
\]
The elements $w_1,w_2,\dots$ of $\Lambda$
determined by the relations
\bel{equ-witt}
\psi_n = \sum_{d | n} dw_d^{n/d} \quad \text{for all } n\in\nn
\end{equation}
also form a free generating set.  This is easy to check using
the following identity:
\[
\sum_{n\geq 0} (-1)^n\lambda_n t^n = \prod_{i\geq 1} (1-x_it)
 = \exp\Big(-\sum_{n\geq 1} \frac{1}{n}\psi_n t^n\Big) = \prod_{n\geq 1}(1-w_nt^n).
\]
The $w_i$ are responsible for the Witt components, as we will see in the
next section.
\end{rmk}

\subsection{}
{\em Remark.}  There is also a description of $\Lambda$ in terms of
the representations of the symmetric
groups~\cite[I.7]{MacDonald:SF}. Let $R_n$ denote the representation
ring of $S_n$, the symmetric group on $n$ letters.  The maps
$S_n\times S_m\to S_{n+m}$, $S_n\to S_n\times S_n$, and $S_n\wr
S_m=S_n\ltimes S_m^n\to S_{mn}$ induce maps between the $R_n$ by
restriction and induction, and these make up a plethory structure on
$\bigoplus_{n\geq0} R_n$ agreeing with that on $\Lambda$.  This is
one natural way to view $\Lambda$ when studying
its action on Grothendieck groups.  (See,
e.g.,~\cite{Deligne:Ex-Ser}).

We do not yet know if similar constructions
in other areas of representation theory also yield plethories.

\subsection{} \lb{sbsc-lambda-p}
{\em $p$-typical symmetric functions.}  Let $p$ be a prime number, and
set $F=\psi_p$.  Then $\zz\ab{F}:=\zz[e,F,F\circ F,\dots]$ is a
subring of $\Lambda$, and because $F$ is ring-like, it is actually a
sub-$\zz$-plethory.  It is also the free plethory on the bialgebra
associated to the monoid $\nn$.  We will denote it $\Psi_p$, and we
will see later that it accounts for the ghost components of the
$p$-typical Witt vectors.

Now let $\Lambda_p$ be the subring of $\Lambda$ consisting of elements
$f$ for which there exists an $i\in\nn$ such that $p^if\in\Psi_p$.
Then $\Lambda_p$ is a sub-$\zz$-plethory of $\Lambda$, and is what we
call the plethory of {\em $p$-typical symmetric functions}.

For all $n\in\nn$, let $\theta_n=w_{p^n}$.
Then (\ref{equ-witt}) becomes
\bel{equ-theta}
F^{\circ n} = \theta_0^{p^n} + \cdots + p^n\theta_n,
\end{equation}
and therefore $\theta_0, \theta_1,\dots$ lie in $\Lambda_p$.
Conversely, because we have
\[
\Lambda = \zz[\theta_0,\theta_1,\dots][w_n \,|\, n\text{ is not a power of }p],
\]
we see $\Lambda_p = \zz[\theta_0,\theta_1,\dots]$.

\subsection{} 
{\em Binomial plethory.}  
Because $\Lambda$ is a $\zz$-plethory, the ring $\zz$ of integers is a $\Lambda$-ring.  The ideal in $\Lambda$ of elements that act as the
constant function $0$ is generated by the set $\{\psi_n-e\,|\, n\geq 1\}$.
The quotient ring is still a plethory, and an action of it on a ring
$R$ is the same as giving $R$ the structure of a $\Lambda$-ring whose
Adams operations are the identity.  This has been shown by Jesse
Elliott (unpublished) to be the same as a binomial $\lambda$-ring 
structure~\cite[p. 9]{Knutson:Lambda} on $R$.

This plethory can also be interpreted as
the set of functions $\zz\to\zz$ that can be expressed as polynomials
with rational coefficients~\cite{Bergman}.

\section{Examples of Witt rings}
\lb{sec-Witt}

Let $k$ be a ring.  Recall that if $P$ is a
$k$-plethory, then $W_P(R)$ denotes the $P$-ring $\Alg{k}(P,R)$.
Because $W_P$ is the right adjoint of the forgetful functor from
$P$-rings to rings, there is a natural map $W_P\to W_P(W_P(R))$, which
in the case of the classical plethories is sometimes called the
Artin--Hasse map.

\subsection{} 
{\em Bialgebras.}
Let $P$ be the free $k$-plethory~\ref{sbsc-free-on-bialg} on a cocommutative
$k$-bialgebra $A$.  Then we have $W_P(B) = \Mod_k(A,B)$.  If $A$ is
finitely generated as a $k$-module, $W_P(B)$ is just $B\tn_k A^*$,
where $A^*$ denotes the dual bialgebra $\Mod_k(A,k)$.  We leave it to
the reader to verify that, in this case,
the map $W_P(B) \longmap W_P(W_P(B))$ is nothing but the
comultiplication map on this bialgebra.  For example, if $A$ is the
group algebra of a finite group $G$, then we have $W_P(B)=B^G$ and the map
above is the map $B^G\to B^{G\times G}=B^G\tn_B B^G$ induced
by the multiplication on $G$.

\subsection{} 
{\em Symmetric functions.}
Because $\Lambda=\zz[\lambda_1,\dots]$, the set $W_\Lambda(B)$ is just
$\prod_{n>0} B$, and it is easy to check that, as a group, we have
$W_\Lambda(B)=1+xB\pser{x}$, where the group operation on the right is
multiplication of power series.  It is also true that if
$1+xB\pser{x}$ is given a $\Lambda$-ring structure as in~\cite[1.1]{AT:Lambda},
then the identification above is an isomorphism of $\Lambda$-rings,
i.e., $W_\Lambda(B)$ is the $\Lambda$-ring of ``big'' Witt vectors.
The proof of this is very straightforward but involves, of course, the
somewhat unpleasant definition of the $\Lambda$-ring structure on
$1+xB\pser{x}$.  Because the whole point of this paper is to move away
from such things, we will leave the argument to the
reader.  The generating set $\{w_1,w_2,\dots\}$
of~\ref{rmk-lambda-basis} allows us to view an element of
$W_\Lambda(B)$ as a (``big'') Witt vector in the traditional 
sense~\cite[17.1.15]{Hazewinkel:Book}.
Under this identification, the map $W_{\Lambda}(B) \longmap
W_{\Lambda}(W_{\Lambda}(B))$ agrees with the usual Artin--Hasse 
map~\cite[17.6]{Hazewinkel:Book}.

If $\Psi$ denotes the sub-plethory $\zz[\psi_n\,|\, n\geq 1]$
of $\Lambda$, then $W_{\Psi}(B)$ is just $\prod_{n>0} B$ as a ring, and
under this identification, the map $W_{\Lambda}(B)\to W_{\Psi}(B)$
is the ghost-component map.

Some early references to the big Witt vectors are
Cartier~\cite{Cartier:Witt-gen} and 
Witt (\cite{Lang:Algebra-I} or \cite[pp.157--163]{Witt:Works}).

\subsection{} 
{\em $p$-typical symmetric functions.}
Because $\Lambda_p=\zz[\theta_0,\dots]$,
the set $W_{\Lambda_p}(B)=\Ring_{\zz}(\Lambda_p,B)$ is naturally
bijective with $B^\nn$.  If we view $B^\nn$ as the set underlying the
ring of $p$-typical Witt
vectors~\cite{Witt:Vectors}\cite[17.1.15]{Hazewinkel:Book}, 
then this bijection is
an isomorphism of rings.  One can write down the
the corresponding $\Lambda_p$-action on $B^{\nn}$, and
we recover the $p$-typical Artin--Hasse map as we did above.
Also as above, if $\Psi_p$ denotes the plethory
$\zz[\psi_p^{\circ\nn}]$, then the natural map $W_{\Lambda_p}(B)\to
W_{\Psi_p}(B)$ is the $p$-typical ghost-component map.

The Teichm\"uller lift can be constructed by considering the monoid
algebra $\zz B$ on the multiplicative monoid underlying $B$.  The ring
$\zz B$ has no additive $p$-torsion, and the map $F\:[b]\mapsto
[b^p]=[b]^p$ ($[\vbl]$ denoting the multiplicative map $B\to \zz B$)
reduces to the Frobenius map modulo $p$.  The ring $\zz[B]$ therefore
(\ref{lemma-cdd}) admits a unique $\Lambda_p$-ring structure where $F$
is the above map.  The canonical ring map $\zz B\to B$ then induces by
adjointness a map $\zz B\to W_{\Lambda_p}(B)$.  In the standard
description, it is $[b]\mapsto (b,0,0,\dots)$, which is of course the
Teichm\"uller lift of $b$.

The following lemma implies that a $\Lambda_p$-ring is the same as
what Joyal calls a $\delta$-ring.  (A comonadic version of this
statement is stated quite clearly in Joyal~\cite{Joyal:Witt}; we
include it only because we will use it later.)

\begin{lemma}
\lb{lemma-cdd}
The $R$ be a $p$-torsion-free ring.  Given an action of $\Lambda_p$ on
$R$, the element $F$ gives an endomorphism of $R$ such that
$F(x)\equiv x^p \mod pR$.  This is a bijection from the set of
actions of $\Lambda_p$ on $R$ to the set of lifts of the Frobenius
endomorphism of $R/pR$.
\end{lemma}
\begin{proof}
Because $R$ is $p$-torsion-free, (\ref{equ-theta}) implies that any
action of $\Lambda_p$ is determined by the endomorphism $F$, and so we
need only show every Frobenius lift comes from some action
of $\Lambda_p$.  

Given a Frobenius lift $f\:R\to R$, Cartier's Dieudonn\'e--Dwork
lemma~\cite[VII\S 4]{Lazard:CFG}
states there is a ring map $R\to W_{\Lambda_p}(R)$ such that the
composite $R\to W_{\Lambda_p}(R)\to W_{\Psi_p}(R)$ sends $r$ to
$(r,f(r),f(f(r)),\dots)$.  This gives a map
$\Lambda_p\bcp R\to R$; to show it is an action we need only check it
is associative.  Because $R$ is $p$-torsion-free it suffices to check
the induced map of $\Psi_p\bcp R\to R$ is an action.  But the
Dieudonn\'e--Dwork lemma implies this map sends
$F^{\circ i}\lcp r$ to $f^{\circ i}(r)$, which is clearly associative.
\end{proof}

\section{Reconstruction and recognition}
\lb{sec-reconstruction}

In preparation for the reconstruction theorem, we generalize the
notions of biring and plethory from $\Ring_k$ to $\Ring_P$ for
non-trivial plethories $P$.  This gives us $P$-$P'$-birings and
$P$-plethories, which reduce to $k$-$k'$-birings and $k$-plethories
when $P=\unt{k}$ and $P'=\unt{k'}$.

Let $P$ be a $k$-plethory and $P'$ a $k'$-plethory, where $k$ and $k'$
are arbitrary rings.

\subsection{}
{\em Functor $\vbl\bcp_{P'}\vbl\:\br{P,P'}\times\Ring_{P'}\to\Ring_{P}$.}
Take $S\in\br{P,P'}$ and $R\in\Ring_{P'}$.
Then $S\bcp_{P'}R$ is
defined to be the coequalizer of the maps of $P$-rings
\begin{gather*}
S\bcp_{k'}P'\bcp_{k'} R \rightrightarrows S\bcp_{k'} R \\
s \lcp \alpha \lcp r \mapsto (s\circ \alpha) \lcp r \\
s \lcp \alpha \lcp r \mapsto s\lcp (\alpha \circ r). 
\end{gather*}

\begin{lemma}
\lb{lemma-biring-adj}
Let $S$ be a $P$-$P'$-biring.  Then the functor
$S\bcp_{P'}\vbl\:\Ring_{P'}\to\Ring_P$ is the left adjoint of the
functor $\Ring_{P}(S,\vbl)$.
\end{lemma}

We leave the proof to the reader.

\begin{prop}
Let $P\to Q$ be a map of plethories.  Then the restriction functor
$\Ring_{Q}\to\Ring_{P}$ preserves limits and coequalizers and has a
left adjoint (``induction'') $Q\bcp_P\vbl$.  If the map $P\to Q$ is an
isomorphism on scalars, it has a right adjoint (``co-induction'')
$\Ring_{P}(Q,\vbl)$ and preserves all colimits.
\end{prop}

\begin{proof}
Because $Q$ is a $Q$-$P$-biring, $Q\bcp_P\vbl$ is left adjoint
(by~\ref{lemma-biring-adj}) to $\Ring_{Q}(Q,\vbl)$, which is the
forgetful functor $\Ring_Q\to\Ring_{P}$. If $P\to Q$ is a map of
$k$-plethories, $Q$ is a $P$-$Q$-biring, so $\Ring_{P}(Q,\vbl)$ is
right adjoint to $Q\bcp_Q\vbl$, the forgetful functor.  It follows
that the forgetful functor preserves limits and, when the rings of
scalars agree, colimits.  It remains to show it always preserves
coequalizers.

Consider the commutative diagram of forgetful functors
\[
\xymatrix{
\Ring_Q \ar[d]\ar[r] & \Ring_{k_Q} \ar[d] \\
\Ring_P \ar[r] & \Ring_{k_P}.
}
\]
The upper functor preserves colimits, and the right-hand functor
preserves coequalizers.  The lower functor reflects isomorphisms and
preserves colimits.  It then follows that the left-hand functor
preserves coequalizers.  
\end{proof}

\begin{rmk}
If $k_P\to k_Q$ is not an isomorphism, $\ep^+$ will fail to descend.
Thus, $Q$ will not be a $k_P$-$k_Q$-biring, let alone a
$P$-$Q$-biring.
\end{rmk}

\subsection{}
A {\em $P$-plethory} is defined to be
a plethory $Q$ equipped with a map $P\to Q$ of
plethories which is an isomorphism on scalars.  A morphism $Q\to Q'$ of
$P$-plethories is a morphism of plethories commuting with the maps
from $P$.

\begin{prop}
\lb{prop-relative-plethory}
$\vbl\bcp_P\vbl$ makes $\br{P,P}$ into a monoidal category with unit
object $P$. Monoids in this category are the same as $P$-plethories.
An action of such a monoid $Q$ on a $P$-ring is the
same as an action of $Q$ on the underlying $k$-ring such that the
action of $Q$ restricted to $P$ is the given one.
\end{prop}
\begin{proof}
The first statement requires no proof. 
Given a monoid $Q$, the structure maps give map $Q\bcp_k Q\to Q\bcp_P Q\to Q$
and $P\to Q$ making it a $k$-plethory. Conversely, a map $P\to Q$
of $k$-plethories makes $Q$ a $P$-$P$-biring and the associativity
condition $Q\bcp_k Q\bcp_k Q\rightrightarrows Q\bcp_k Q\to Q$
implies that $Q\bcp_k P \bcp_k Q\rightrightarrows Q\bcp_k Q\to Q$
commutes, so composition descends to $Q\bcp_P Q\to Q$.

Similarly, an action of $Q$ on the underlying $k$-ring of a
$P$-ring $A$ is a map $Q\bcp_kA\to A$, and it
descends to a $P$-action $Q\bcp_PA\to A$ because
$Q\bcp_kP\bcp_kA\rightrightarrows Q\bcp_kA\to A$ commutes.
\end{proof}

\subsection{}
\lb{sbsc-recon-setup}
Now let $\mathsf{C}$ be a category that has
all limits and colimits, and let $U\:\mathsf{C}\to \Ring_{P}$ be a
functor that has a left adjoint $F$. We also assume $U$ reflects
isomorphisms, that is, a morphism $f$ is an isomorphism if and only if
$U(f)$ is an isomorphism.  Set $Q=UF(P)$.
Let $U'$ be the composite of $U$ with the forgetful functor from
$\Ring_P$ to the category of sets.

\subsection{} \lb{sbsc-tannaka-const}
{\em $k$-Plethory structure on $Q$ when $U$ has a right adjoint.}
Suppose $U$ has a right adjoint $W$.  The functor $UW$ is represented
by $Q$: $UW(A)=\Ring_P(P,UW(A))=\Ring_{P}(UF(P),A)$, and this gives
$Q$ the structure of a $P$-$P$-biring (\ref{sbsc-P-action}). The
composite $UW$ of adjoints is a comonad, and so its adjoint
$Q\bcp_P\vbl$ is a monad. By~\ref{prop-relative-plethory}, $Q$ is a
$k$-plethory with a map $P\to Q$.

Given an object $A$ of $\mathsf{C}$, the adjunction gives an action of
$UF(\vbl)=Q\bcp\vbl$ on $U(A)$, and hence we have a functor
$\mathsf{C}\to\Ring_Q$ between categories over $\Ring_P$.

\begin{thm}
\lb{thm-pre-tannaka}
If $U$ has a right adjoint $W$, then the functor
$\mathsf{C}\to\Ring_Q$ is an equivalence of categories over
$\Ring_P$.
\end{thm}

\begin{proof}
Beck's theorem~\cite{MacLane:CWM}.
\end{proof}

\subsection{}
\lb{sbsc-U'}
Let $k'$ be the $P$-ring $UF(k)$, and let $P'$ be the $k'$-plethory
$k'\tn_k P$.  Because $F(k)$ is the initial object, $U$ factors as a
functor $U'\:\mathsf{C}\to\Ring_{P'}$ followed by the forgetful
functor $V\:\Ring_{P'}\to\Ring_P$.  The functor
$U'$ has a left adjoint $F'$ given by descent: if $A$ is a $P'$-ring,
then $FV(A)$ has two maps from $F(k')=FUF(k)$, one from applying
$FV$ to the initial map $k'\to A$ and the other given by the
composite
\[
FUF(k)\to F(k) \to FV(A),
\]
where the first map is the adjunction and the second is the initial
map.  Let $F'(A)$ denote the coequalizer of $F(k')\rightrightarrows
FV(A)$.

\begin{thm}
If $P\to Q$ is a map of plethories and $U$ is the forgetful functor
$\Ring_Q\to\Ring_P$, then $U'$ of~\ref{sbsc-U'} has a right adjoint.
Conversely, suppose $U'\:\mathsf{C}\to\Ring_P$ has a right adjoint,
and let $Q$ be the $k'$-plethory $U'F'(P')$
of~\ref{sbsc-tannaka-const}.  Then the functor $\mathsf{C}\to\Ring_Q$
is an equivalence of categories over $\Ring_{P'}$.
\end{thm}

\begin{proof}
Apply~\ref{thm-pre-tannaka} to $U'$.
\end{proof}

\begin{rmk}
In practice, it is quite easy to check the existence of $F$ and $W'$
using Freyd's theorem from category theory.
\end{rmk}

\section{$P$-ideals}

Let $P$ be a $k$-plethory, and let $P_+$ denote the kernel of
$\ep^+\:P\to k$.

\subsection{}
An ideal $I$ in a $P$-ring $R$ is called a {\em (left) $P$-ideal} if
there exists an action of $P$ on $R/I$ such that the map $R\to R/I$ of
rings is a map of $P$-rings.  If such an action exists, it is unique,
and so being a $P$-ideal is a property of, rather than a structure on, a
subset of $R$.

\begin{prop}
\lb{prop-P-ideal}
Let $I$ be an ideal in a $P$-ring $R$.  Then the following are
equivalent:
\begin{enumerate}
\item $I$ is a $P$-ideal
\item $I$ is the kernel of a morphism of $P$-rings 
\item $P_+\circ I \subseteq I$
\item $I$ is generated by a set $X$ such that
$P_+\circ X\subseteq I$.
\end{enumerate}
\end{prop}

The proof is in \ref{pf-P-ideal}.

Given any subset $X$ of $P$, it is therefore reasonable to call the
ideal generated by $P_+\circ X$ the {\em $P$-ideal generated by $X$}.

\subsection{}
\lb{sbsc-binary}
Elements of $P\otimes P$ give binary
operations on any $P$-ring $R$ by $(\alpha\tn\beta)(r,s) =
\alpha(r)\beta(s)$ and extending linearly.

\begin{lemma}
\lb{lemma-blah}
Let $R$ be a $P$-ring, $I$ an ideal in $R$ and $X$ a subset of $R$.
Assume that for all $x\in X$ and $f\in P_+$, we have $f(x)\in I$.
Then for all $t\in P\tn P_+$ and all $(r,i)\in R\times I$, 
we have $t(r,i)\in I$.  
\end{lemma}
\begin{proof}
Since $t\in P\tn P_+$, it may be expressed as $t=\sum t'_j\tn t''_j$
with $t''_j\in P_+$, so that $t''_j$ preserves $I$. 
Then for $(r,i)\in R\times I$, $t(r,i)=\sum t'_j(r)t''_j(i)\in I$.
\end{proof}
Typical applications will use $X=I$, a $P$-ideal.

\begin{lemma}
\lb{lemma-delta*}
Let $S$ be a $k$-$\zz$-biring.  Then
$\Delta^+(S_+)$ is contained in $S_+\tn S+S\tn S_+$, and
$\Delta^{\times}(S_+)$ is contained in $S_+\tn S_+$.
\end{lemma}
\begin{proof}
$S$ is a ring object in the opposite of $\Alg{k}$;
the ring identity $0+0=0$ translates into the identity
$(\ep^+\tn\ep^+)\circ \Delta^+ = \ep^+$, which is clearly equivalent
to the first statement.  The second statement is similarly just a
coalgebraic translation of a ring identity.
Let $W$ denote the ring object corresponding to $S$ in the opposite category.
Then the commutativity of the following two diagrams is equivalent:
\[
\begin{array}{ccc}
\xymatrix{
W & W\times W\ar[l]_-{{\times}} \\
0\ar[u] & W\,{\scriptstyle\amalg}\, W\ar[l]\ar[u]_{\id\times 0 \,{\scriptscriptstyle\amalg}\, 0\times \id}.
}
& \hspace{5mm} &
\xymatrix{
S\ar[r]^{\Delta^{\times}}\ar[d]^-{\ep^+} & S\tn S
  \ar[d]^{\id\tn\ep^{+} \times \ep^+\tn \id} \\
k\ar[r] & S\times S
}
\end{array}
\]
But the commutativity of the first is just a restatement of the ring
identity $0\cdot x = x\cdot 0 = 0$.
We therefore have
$$
\Delta^{\times}(S_+) \subseteq \ker\Big(S\tn S\to S\times S\Big)
 = S_+\tn S_+.
$$
\end{proof}

\subsection{}
\lb{pf-P-ideal}
{\em Proof of \ref{prop-P-ideal}.}
(1)$\Rightarrow$(2) and (3)$\Rightarrow$(4) are clear.

(2)$\Rightarrow$(3):
$P_+$ preserves the set $\{0\}$ in $k$ and, thus, in any $P$-ring; it
therefore must preserve its preimage under a morphism of $P$-rings.

(3)$\Rightarrow$(1):
If $I$ is preserved by $P_+$, we must put a $P$-ring structure on
$R/I$ so that $R\to R/I$ is a morphism of $P$-rings. The action must
be $p(r+I)=p(r)+I$; it is necessary only to check that this is
well-defined. The kernel of $\id_P\tn \ep^+\:P\tn P\to P$ is $P\tn
P_+$, and so by the counit condition, we have
$\Delta^+p- p\tn1\in P\tn P_+$ for all $p\in P$.
For any $i\in I$, we have $p(r+i)-p(r)=(\Delta^+p-p\tn1)(r,i)$.
By~\ref{lemma-blah}, the right-hand side of this equality is in $I$,
and so the action is well-defined.

(4)$\Rightarrow$(3): 
Consider the set $J$ of elements of $I$ that are sent into $I$ by all
elements of $P_+$. If $f\in P_+$, then $\Delta^+f\in P_+\tn P + P\tn
P_+$.  Thus for $j,k\in J$, lemma~\ref{lemma-blah} implies $f(j+k)\in
I$ and hence $j+k\in J$.  Similarly, $\Delta^\times f\in P_+\tn
P_+\subset P\tn P_+$, and so for $r\in R$ and $j\in J$, we have
$f(rj)\in I$ and hence $rj\in J$.  Therefore $J$ is an ideal, and if
a generating set for $I$ is sent by $P_+$ into $I$, we have $I=J$.
So all of $I$ is preserved by $P_+$.
\qed

\begin{prop}
\lb{prop-ideal-product}
Let $I$ and $J$ be $P$-ideals in a $P$-ring $A$.  Then $IJ$ is a $P$-ideal.
\end{prop}
\begin{proof}
It is sufficient to check $f(xy)\in IJ$ for all $f\in P_+$, $x\in I$,
and $y\in J$ because such $xy$ form a generating set. We can write
$\Delta^\times f=\sum f^{[1]}_i\tn f^{[2]}_i$ with 
$f^{[1]}_i,f^{[1]}_i\in P_+$, and so we have
$f(xy)=\sum f^{[1]}_i(x)f^{[2]}_i(y)\in IJ$.
\end{proof}

\section{Two-sided ideals}

Let $P$ be a $k$-plethory, and let $P'$ be a $k'$-plethory.

\subsection{}
\lb{sbsc-bi-ideal-def}
An ideal $J$ in a $k$-$k'$-biring $S$ is called a {\em $k$-$k'$-ideal}
if the quotient $k$-ring $S/J$ admits the structure of a
$k$-$k'$-biring.  This is clearly equivalent to $S/J$ being, in the
opposite of $\Alg{k}$, a sub-$k'$-ring object of $S$, and so if $S/J$
admits such a structure, it is unique.  This is also
equivalent to the existence of a generating set $X$ of $J$ 
such that, in the notation of~\ref{sbsc-biring}, we have
\begin{enumerate}
\item $\Delta^+_S(X) \subseteq S\tn J + J\tn S$,
\item $\Delta_S^\times(X) \subseteq S\tn J + J\tn S$, and
\item $\beta_S(c)(X)=0$ for all $c\in k'$.
\end{enumerate}

\subsection{}
A $k$-$k'$-ideal $J$ in a $P$-$P'$-biring $S$ is called a {\em $P$-$P'$-ideal}
if there exists a $P$-$P'$-biring
structure on the quotient $k$-$k'$-biring $S/J$ such that $S\to S/J$
is a map of $P$-$P'$-birings.  If such an action exists, it is unique,
and so as was the case for $P$-ideals, being a $P$-$P'$-ideal is a
property, rather than a structure.

\begin{prop}
Let $J$ be a $k$-$k'$-ideal in a $k$-$k'$-biring $S$.  Then the following are
equivalent
\begin{enumerate}
\item $J$ is a $P$-$P'$-ideal
\item $J$ is the kernel of a map of $P$-$P'$-birings
\item $P_+\circ J\circ P' \subseteq J$
\item $J$ is generated by a set $X$
such that $P_+\circ X\circ P'\subseteq J$.
\end{enumerate}
\end{prop}

The asymmetry in (3) is due to the traditional definition of ideal.
If we took a more categorical approach and considered, instead of kernels
of maps $R\to S$ of $k$-rings, the fiber products $R\times_S k$, the
$P_+$ in (3) would become a $P$.

\begin{proof}
As in~\ref{prop-P-ideal}, the only implication that requires proof is
$(4)\Rightarrow(1)$. 

So, assume (4).  By~\ref{prop-P-ideal}, $J$ is a $P$-ideal; and by
assumption, $J$ is a $k$-$k'$-ideal.  Therefore $S/J$ is a
$P$-$k'$-biring.  For all $s\in S$,$j\in J$,$ f\in P'$, we have
\[
(s+j)\circ f = s\circ f + j\circ f \equiv s\circ f \mod J,
\]
and so the right $P'$-action descends to $S/J$.
\end{proof}

\subsection{}
\lb{sbsc-plethory-quot}
If $J$ is a $P$-$P$-ideal in $P$ itself, then this proposition implies
$P/J$ is a $P$-plethory in the sense that the $P$-$P$-biring structure
on $P/J$ extends to a unique $P$-plethory structure on $P/J$.

\begin{prop}
\lb{prop-biring-colimits}
The category $\br{P,P'}$ of $P$-$P'$-birings has all colimits, and the
forgetful functor $\br{P,P'}\to\Ring_P$ preserves them.
\end{prop}
\begin{proof}
Given a diagram $C$ of $P$-$P'$-birings, its colimit $S$ in the
category of $P$-rings has the property that for any $P$-ring $R$, the
set $\Ring_P(S,R)$ is the limit of the sets $\Ring_P(T_c,R)$, where
$c$ ranges over $C$.  Because each $\Ring_P(T_c,R)$ is a $P'$-ring and
the maps are $P'$-equivariant, $\Ring_P(S,R)$ is a $P'$-ring.  Thus,
by a remark in~\ref{sbsc-P-action}, $S$ has a unique $P$-$P'$-biring
structure making the maps $T_c\to S$ maps of $P$-$P'$-birings, which
was to be proved.
\end{proof}

\subsection{} \lb{sbsc-free-on-pointed-biring}
{\em Free plethory on a pointed biring.}
The free $P$-plethory $Q$ on a $P$-$P$-biring $S$ can be constructed as
in~\ref{sbsc-free-on-biring}.
It comes equipped with a map $P\to Q$ of $k$-plethories.

Now let $f\:P\to S$ be a map of $P$-$P$-birings. (This is equivalent
to specifying an element $s_0\in S$ such that $p\circ s_0=s_0\circ p$
for all $p\in P$.)  Then the free plethory on the pointed biring
$S$ is the coequalizer (\ref{prop-biring-colimits}) of the two
$Q$-$Q$-biring maps $Q\bcp P\bcp Q\rightrightarrows Q$ induced by
sending $e\lcp\alpha\lcp e$, on the one hand, to $\alpha\in P=S^{\lcp
0}$ and, on the other, to $f(\alpha)\in S^{\lcp 1}$.
By~\ref{sbsc-plethory-quot}, $Q$ is a $k$-plethory.  It is the initial
object among $P$-plethories $P'$ equipped with a map $S\to P'$ such
that the composite $P\to S\to P'$ agrees with the structure map $P\to
P'$.  An action of this plethory on a $k$-ring $R$ is the same as
an action of $P$ on $R$ together with a map $S\bcp R\to R$ such that
$f(p)\lcp r\mapsto p(r)$ for all $p\in P, r\in R$.

At this point, it is quite easy to give an explicit construction of $\Lambda_p$
that does not rely on symmetric functions.  Let
$S=\zz[e,\theta_1]$ be the $\unt{\zz}$-pointed $\zz$-$\zz$-biring
determined by
\begin{align}
\lb{equ-theta-coprod}
\Delta^+\: \theta_1 &\mapsto \theta_1\otimes 1 + 1\otimes\theta_1 
            - \sum_{i=1}^{p-1}\frac{1}{p}\binom{p}{i}e^i\otimes e^{p-i} \\
\Delta^\times\: \theta_1 &\mapsto e^p\otimes \theta_1 + \theta_1\otimes
            e^p + p\theta_1\otimes\theta_1 
\end{align}
Then Cartier's Dieudonn\'e--Dwork lemma implies $\Lambda_p$ is the
free $\zz$-plethory on $S$.  Of course, this is just a plethystic
description of Joyal's approach~\cite{Joyal:Witt} to the $p$-typical
Witt vectors.

\subsection{}
\lb{sbsc-mod-pointed}
The following asymmetric variant of this construction will be used in
section~\ref{sec-amplification}.  Let $P_0$ be a $k$-plethory, let $P$ be a
$P_0$-plethory, let $S$ be a $P_0$-$P$-biring, and let $g\:P\to S$ be a
map of $P_0$-$P$-birings.
Let $Q$ denote the free $P_0$-plethory on $S$ viewed as a pointed $P_0$-$P_0$-biring.
Then we have two maps of $P_0$-$P_0$-birings $S\bcp_{P_0}
P\rightrightarrows Q$ given by $s\lcp \alpha\mapsto s\bcp
g(\alpha)\in S^{\lcp 2}$
and $s\lcp \alpha \mapsto s\circ \alpha\in S^{\lcp 1}$.  These then
induce two maps of $Q$-$Q$-birings $Q\bcp_{P_0} S\bcp_{P_0}
P \bcp_{P_0} Q \rightrightarrows Q$.
The coequalizer $T$ of these maps is a $P_0$-plethory
(\ref{sbsc-plethory-quot}), but the two maps $P\to Q$ become equal
in $T$, and so $T$ is in fact a $P$-plethory.  An action of $T$ on a
ring $R$ is the same as an action of $P$ on $R$ together with a map
$S\bcp_{P} R\to R$ such that $g(\alpha)\lcp r \mapsto \alpha\circ r$.

\section{Amplifications over curves}
\lb{sec-amplification}

Let $\sO$ be a Dedekind domain, and let $\m$ be an ideal; let $k$
denote the residue ring $\sO/\m$, and let $K$ denote the subring of
the field of fractions of $\sO$ consisting of elements that are integral at
all maximal ideals not dividing $\m$.  The $\m$-torsion submodule of
an $\sO$-module $M$ is the set of $m\in M$ for which there exists an
$n\in\nn$ such that $\m^n m=0$.  We say an $\sO$-module is
$\m$-torsion-free if its $\m$-torsion submodule is trivial, or equivalently, if
it is flat locally at each maximal ideal dividing $\m$.

Now let $P$ be an $\sO$-plethory that is $\m$-torsion-free,
let $Q$ be a $k$-plethory, and let $f\:P\to Q$ be a
surjective map of plethories agreeing with the canonical map on
scalars.  A {\em $P$-deformation of a $Q$-ring} is an
$\m$-torsion-free $P$-ring $R$ such that the action of $P$ on
$k\tn R$ factors through an action of $Q$ on $k\tn R$. (Note that
because $P\to Q$ is surjective, it can factor in at most one way.)
A morphism of $P$-deformations of $Q$-rings
is by definition a morphism of the underlying $P$-rings.

The purpose of this section is then to construct an $\sO$-plethory
$P'$, {\em the amplification of $P$ along $Q$}, such that
$\m$-torsion-free $P'$-rings are the same as $P$-deformations of
$Q$-rings.  It is constructed simply by adjoining $\m^{-1}\tn I$
to $P$, where $I$ is the kernel of the map $P\to Q$, and so it is
analogous to an affine blow-up of rings.  Note however that there are
some minor subtleties involved in adjoining these elements because a
plethory involves co-operations, not just operations, and because we need to
know how to compose elements of $P$ with elements of $\m^{-1}\tn I$, but
$P$ may not even act on $K$, let alone preserve $\m$.  

\begin{thm}
\lb{thm-amplification}
The $P$-plethory $P'$ of~\ref{sbsc-amplification} is $\m$-torsion-free,
and the forgetful functor from the full category of $\m$-torsion-free
$P'$-rings to $\Ring_P$ identifies it with the category of $P$-ring
deformations of $Q$-rings.  Furthermore, $P'$ has the following
universal property: Let $P''$ be a $P$-plethory whose underlying
$P$-ring is a $P$-deformation of a $Q$-ring.  Then there is a unique
map $P'\to P''$ of $P$-rings commuting with the maps from $P$, and
this map is a map of $P$-plethories.
\end{thm}

\begin{cor}
\lb{cor-amplification}
Let $P''$ be a $P$-plethory with the property that the forgetful
functor from the full category of $\m$-torsion-free $P''$-rings to
$\Ring_P$ identifies it with the category of $P$-ring deformations of
$Q$-rings.  Then there is a unique map $P''\to P'$ of $P$-rings; this
map is a map of $P$-plethories, and it identifies $P'$ with the
largest $\m$-torsion-free $P''$-ring quotient of $P''$.
\end{cor}

We prove these at the end of this section.  Note that either the
theorem or the construction of~\ref{sbsc-amplification} implies
amplification is functorial in $P$ and $Q$.

\subsection{} \lb{rmk-mod}
{\em Remark.}  
As always, either universal property determines $P'$ uniquely up to
unique isomorphism.  The final statement of the corollary determines
it without any mention of universal properties: it is the unique
$\m$-torsion-free $P$-plethory such that the forgetful functor
identifies $\m$-torsion-free $P'$-rings with $P$-deformations of
$Q$-rings.

One could also describe the category of all $P'$-rings as the category
obtained from the category of $P$-deformations of $Q$-rings (i.e.,
$\m$-torsion-free $P'$-rings) by adjoining certain colimits.  This would give
another satisfactory approach to the functor of $p$-typical Witt vectors
circumventing any discussion of plethories.

\begin{lemma}
\lb{lemma-flat-quot}
Let $T$ be an $\sO$-plethory.
Then the $T$-ideal in $T$ generated by the $\m$-torsion ideal 
is a $T$-$T$-ideal.
\end{lemma}
\begin{proof}
Let $I$ denote the ideal of $\m$-torsion in $T$, and let $J$ denote
the $T$-ideal it generates.  First we show $I$ is an
$\sO$-$\sO$-ideal.  Because $I$ is $\m$-torsion, the ideal $T\tn I +
I\tn T$ is contained in the $\m$-torsion ideal of $T\tn T$.
But this containment is actually an equality: because $T/I$ is
$\m$-torsion-free and because $\sO$ is a Dedekind domain, $T/I\tn
T/I$ is $\m$-torsion-free.
It therefore follows that $\Delta^+(I)$ and
$\Delta^{\times}(I)$ are both contained in $T\tn I + I\tn T$.
And last, $\beta(c)(I)$ is zero because it is torsion but $\sO$ is
torsion-free.  By~\ref{sbsc-bi-ideal-def}, the ideal $I$ is an
$\sO$-$\sO$-ideal.

Now we show $J$ is a $T$-$T$-ideal.  It is a $T$-ideal by definition,
and so we need only show $J \circ T\subseteq J$, and in fact only
$I\circ T\subseteq J$.  So take $i\in I$ and $\alpha\in T$.  Then
there is some $n\in\nn$ such that $\m^n i=0$, and for every $x\in\m^n$,
we have $x(i\circ\alpha) = (xi)\circ\alpha =0$.
\end{proof}

\subsection{} \lb{sbsc-flat-quot}
{\em Maximal $\m$-torsion-free quotient of an $\sO$-plethory.}  Let $T_0$ be
an $\sO$-plethory, let $J$ denote the $T_0$-ideal generated by the
$\sO$-torsion.  By~\ref{lemma-flat-quot} and~\ref{sbsc-plethory-quot},
the quotient $T_1=T/J$ is an $\sO$-plethory.  Let $T_2$ be the same
construction applied to $T_1$, and so on.  Then the colimit of the sequence
\[
T_0\longmap T_1 \longmap \cdots
\]
in the category of $T_0$-$T_0$-birings
(\ref{prop-biring-colimits}) is clearly the largest $\m$-torsion-free
$T_0$-ring quotient of $T_0$.  It is an $\sO$-plethory because
it is a quotient $T_0$-$T_0$-biring of $T_0$.

Note that $\m$-torsion-free $T_0$-rings are the same as
$\m$-torsion-free $T'$-rings.

\subsection{} 
\lb{sbsc-amplification}
{\em Amplification $P'$ of $P$ along $Q$.}  Let $I$ denote the kernel
of the map $P\to Q$, and let $S$ denote the sub-$\sO$-ring of $K\tn
P$ generated by $\m^{-1}\tn I$. (Here, all tensor products are over
$\sO$, and as usual $\m^{-1}$ denotes the
$\sO$-dual of $\m$ viewed as a submodule of $K$.)  Note that we have
$1\otimes P\subseteq S$ and also that $K\tn P$ is a $\unt{K}$-$P$-biring,
but it need not be a $K$-plethory.

The $K$-$\sO$-biring structure on $K\tn P$ induces an
$\sO$-$\sO$-biring structure on $S$ as follows: Let $\Delta$ denote
either $\Delta^+$ or $\Delta^{\times}$, and let $\Delta_K$ denote
$\id_K\tn\Delta$.  Then we have
\begin{align*}
\Delta_K(\m^{-1}\tn I) 
   &\subseteq \m^{-1}\tn\Delta(I) \\
   &\subseteq \m^{-1}\tn (P\tn I + I\tn P).
\end{align*}
Identifying $K\tn P\tn P$ with $(K\tn P)\tn (K\tn P)$, we have
\[
\Delta_K(\m^{-1}\tn I) 
   \subseteq (1\tn P)\tn (\m^{-1}\tn I)+(\m^{-1}\tn I)\tn (1\tn P)
   \subseteq S\tn S.
\]
Because $\Delta_K$ is an $\sO$-ring map, it follows that
$\Delta(S)\subseteq S\tn S$.
Similarly, if $\ep$ denotes either the additive or
multiplicative counit and $\ep_K=\id_K\tn\ep$ , then
\[
\ep_K(\m^{-1}\tn I) = \m^{-1}\tn \ep(I)
 \subseteq \m^{-1}\tn \m
 = \sO,
\]
and as above, we have $\ep(S)\subseteq \sO$.
The properties necessary for this data to give a $\sO$-$\sO$-biring
structure on $S$ follow from the $K$-$\sO$-biring properties on $K\tn P$.

Because $I$ is preserved by the right action of $P$, so is $S$, and
therefore $S$ has a $\unt{\sO}$-$P$-biring structure.  Let
$T$ be the construction of~\ref{sbsc-mod-pointed} applied to the
$\sO$-plethory $P$, the $\unt{\sO}$-$P$-biring $S$, and the inclusion
map $P\to S$.

Finally, let $P'$ denote the maximal $\m$-torsion-free quotient of $T$ 
(\ref{sbsc-flat-quot}).  It is a $P$-plethory because $T$ is.

\begin{lemma}
\lb{lemma-amplification}
Let $R$ be an $\m$-torsion-free $P$-ring.  Then the action of $P$ on
$R$ factors through at most one action of $P'$, and one
exists if and only if $R$ is a $P$-deformation of a $Q$-ring.
\end{lemma}

\begin{proof}
Suppose the action of $P$ on $R$ prolongs to two actions $\circ_1$ and
$\circ_2$ of $P'$.  For any $\alpha\in P'$ and $r\in R$, we want to
show $\alpha\circ_1 r = \alpha\circ_2 r$.  Because $T$ surjects onto
$P'$, it is enough to show this for $\alpha$ in $T$ and, because $S$
generates $T$, even in $S$.  But $S$ is a subset of $K\tn P$; so take
some $n\in\nn$ such that $\m^n\alpha\subseteq P$.  Then
\[
x(\alpha\circ_1 r) = (x\alpha)\circ_1 r = (x\alpha)\circ_2 r =
x(\alpha\circ_2 r)
\]
for all $x\in\m^n$. But because $R$ is
$\m$-torsion-free, we have $\alpha\circ_1 r = \alpha\circ_2 r$, and so
there is at most one compatible action of $P'$ on $R$.

The action of $P$ on $R/\m R$ factors through $Q$ if and only if
$I\circ R\subseteq \m R$.  This is equivalent to
$(\m^{-1}\tn I)\circ R\subseteq R$ under the map
\[
(K\tn P)\bcp R = K\tn (P\bcp R) \longlabelmap{\circ} K\tn R,
\]
which is in turn equivalent to $S\circ R\subset R$.
Because $R$ is $\m$-torsion-free and because $K\otimes S = K\otimes
P$, this is then equivalent to the
existence of some map $\circ'\:S\bcp_P R\to R$ of $\sO$-rings such
that $p \circ' r = p\circ r$ for all $p\in P, r\in R$.  By
\ref{sbsc-mod-pointed}, this is equivalent to an action of $T$ on
$R$ that is compatible with the given action of $P$, and this is
equivalent to such an action of $P'$ on $R$.
\end{proof}

\subsection{}
{\em Proof of~\ref{thm-amplification}.} 
$P'$ is $\m$-torsion-free by construction.

The forgetful functor is clearly faithful, and
lemma~\ref{lemma-amplification} implies its image is as stated.  To see
it is full, let $R$ and $R'$ be $\m$-torsion-free $P'$-rings and let
$f\:R\to R'$ be a map of $P$-rings.  We need to check $f(\alpha\circ
r) = \alpha\circ f(r)$ for all $\alpha\in P'$ and $r\in R$.  As in the
proof of~\ref{lemma-amplification}, it is enough to show this for
$\alpha$ in $S$, where the equality follows because $R'$ is
$\m$-torsion-free.  This proves the functor is fully faithful.

Let $P''$ be as in the universal property.  By the previous paragraph,
the action of $P$ on $P''$ extends uniquely to an action of $P'$; and
because $P'$ is the free $P'$-ring on one element, there is a unique
map of $P'$-rings $P'\to P''$ sending $e$ to $e$.  Again by the
previous paragraph, we see there is a unique map $P'\to P''$ of
$P$-rings sending $e$ to $e$, that is, commuting with the maps from
$P$.

To show this is a map of $P$-plethories, it is enough to show there
exists some map $P'\to P''$ of $P$-plethories.  Because $P''$ is
$\m$-torsion-free, such a map is the same as a map $T\to P''$ of
$P$-plethories, and this is the same as
a map $S\to P''$ of $\unt{\sO}$-$P$-birings respecting the maps from
$P$.  Because $P''$ is $\m$-torsion-free, there is at most one such
map, and there is exactly one if the map $P\to P''$ sends $I$ to $\m
P''$.  But this is just another way of saying the $P$-ring underlying
$P''$ is a $P$-deformation of a $Q$-ring, and that fact we are given.
\qed

\subsection{}
{\em Proof of~\ref{cor-amplification}.}  
Replacing $P''$ with its maximal $\m$-torsion-free quotient
(\ref{sbsc-flat-quot}), we can assume $P''$ is $\m$-torsion-free.
Then $P''$ and $P'$ are both initial objects in the category of
$P$-deformations of $Q$-rings and so are uniquely isomorphic.  The
universal property of the theorem applied to $P''$ then implies this
isomorphism is a map of $\sO$-plethories.
\qed

\subsection{}
\lb{prop-amplification-K}
Suppose $K$ admits a $P$-action.  Then $K$ is trivially a
$P$-deformation of a $Q$-ring and,
by~\ref{thm-amplification}, has a unique compatible $P'$-action.
By~\ref{sbsc-base-change}, there is a canonical
$K$-plethory structure on $K\otimes P'$.

\begin{prop*}
If $K$ admits a $P$-action, the map $K\tn P\to K\tn P'$ is an
isomorphism of $K$-plethories.  Moreover, under
this identification, $P'$ is the $\sO$-subring of $K\tn P$ generated
by the $\circ$-words in the elements of $\m^{-1}\tn I$.
\end{prop*}

\begin{proof}
To show the first statement, it is enough to show the map induces an
equivalence between $\Ring_{K\tn P}$ and $\Ring_{K\tn P'}$.  But a
$K\tn P$-ring structure on a $K$-ring $R$ is the same
(by~\ref{sbsc-base-change}) as an action of $P$ on $R$, and because
$R$ is trivially a $P$-deformation of a $Q$-ring, this is the same as
a $P'$-action, which (by~\ref{sbsc-base-change} again) is the same as a $K\tn
P'$-ring structure on $R$.

Because $P'$ is $\m$-torsion-free, it is naturally an $\sO$-subring of
$K\tn P'=K\tn P$.
Since $P'$ is the surjective image of the free plethory on the biring
$S$, it is the smallest $\sO$-ring in $K\tn P$ containing $\m^{-1}\tn I$
(and hence $S$) and closed under composition.
\end{proof}

\section{The cotangent algebra}
\lb{sec-cotan1}

By the structure of an {\em algebra over $k$} on an $\zz$-algebra $A$,
we mean simply a morphism $k\to A$ of $\zz$-algebras.  The
image need not be central.  These form a category in the obvious way.

For any $k$-$k'$-biring $S$, write $\ct{S}$ for the $k$-module
$S_+/S_+^2=\ker(\ep_S^+)/\ker(\ep_S^+)^2$.  It is called the cotangent
space of $S$.  The purpose of this
section is to show that the cotangent space is naturally a
$k$-$k'$-bimodule and, especially, the cotangent space of a $k$-plethory is
naturally an algebra over $k$.
We do this by showing that if $S'$ is a
$k'$-$k''$-biring, then $\ct{S\bcp S'}=\ct{S}\tn_{k'}\ct{S'}$.  Thus,
when $k=k'=k''$, the cotangent space is a monoidal functor, so it
sends plethories (monoids in the category of $k$-$k$-birings) to
algebras over $k$ (monoids in the category of $k$-$k$-bimodules).

First we show all elements of $S_+$ are additive up to second order:

\begin{lemma}
\lb{lemma-add-1st-order}
Let $J$ denote the kernel of the map
$\ep^+\tn\ep^+\: S\tn_k S\to k$.  Then
for all $s\in S_+$, we have
$\Delta^+(s) \equiv s\tn 1 + 1\tn s \mod J^2$.
\end{lemma}

\begin{proof}
The cotangent space functor takes coproducts in $\Ring_k$ to
coproducts of $k$-modules and (hence) takes cogroup objects to cogroup
objects.  In particular, we have an identification
$J/J^2=\ct{S}\oplus\ct{S}$, and under this identification, the map
$\ct{S}\to\ct{S}\oplus\ct{S}$ of cotangent spaces induced by
$\Delta^+$ makes $\ct{S}$ a cogroup in the category of $k$-modules.
But the only cogroup structure on a $k$-module is the diagonal map.
\end{proof}

\begin{prop}
Consider the right action of $k'$, as a monoid, on $S$ given by
setting $s\cdot c$ to be the image of $s\lcp ce$ under the
identification $S\bcp_{k'} \unt{k'}=S$.  (Explicitly, $s\cdot c=\sum
\beta(c)(s^{[1]}_i)s^{[2]}_i$.)  Then this action preserves $S_+$ and
descends to $\ct{S}$, and the resulting action makes the $k$-module
$\ct{S}$ a $k$-$k'$-bimodule.
\end{prop}
\begin{proof}
The action preserves $S_+$ since $\ep^+(s)=s\cdot 0$.
Because it acts by ring endomorphisms, it also preserves $S_+^2$, and
thus it descends to $\ct{S}$.
By~\ref{lemma-add-1st-order}, $k'$ acts not just as a monoid, but
as a ring.  It commutes with the $k$-action because for any $b\in k$, we have 
$(bs)\lcp (ce) = (b\lcp (ce))(s\lcp(ce)) = b(s\lcp(ce))$ in $S\bcp\unt{k'}$.
\end{proof}

\begin{prop}
The map $k\to\ct{\unt{k}}$ given by $c\mapsto ce$ is an isomorphism of
$k$-$k$-bimodules.  If $S$ is a $k$-$k'$-biring and $S'$ a
$k'$-$k''$-biring, then the map
$\ct{S}\otimes_{k'}\ct{S'}\to\ct{S\bcp_{k'}S'}$ given by
$s\otimes s'\mapsto s\lcp s'$ is well-defined and an isomorphism of
$k$-$k'$-bimodules.
\end{prop}
\begin{proof}
The first statement follows immediately from the definition
(\ref{sbsc-ident}) of $\unt{k}=k[e]$.

Now we will show the second map is well-defined.
Note that $\ep^+(s\lcp s')=s(s'(0))$, where $\alpha(c)$ denotes
$\beta(c)(\alpha)$.  Thus if $s\in S_+$ and $s'\in S'_+$, then $s\lcp
s'\in (S\bcp S')_+$, and so we have a well-defined map $S_+\times
S'_+\to \ct{S\bcp S'}$. This map is clearly additive in the first variable
and is additive in the second by~\ref{lemma-add-1st-order}. Thus to
check that it descends to $\ct{S}\times\ct{S'}$, we need only show
$s\lcp s'\in (S\bcp S')_+^2$ for $s\in S_+^2$ and $s\lcp s'\in
(S\bcp S')_+^2$ for $s'\in (S')_+^2$.  The first is clear, for ring
operations come out of the left side of the composition product. For
the second, $s'$ may be a sum of products, but up to second order,
sums also come out of the right side (by~\ref{lemma-add-1st-order}),
and so we may assume $s'=s'_1s'_2$, $s'_i\in S'_+$. Then $s\lcp
s'_1s'_2=\Delta^{\times}s(s'_1,s'_2)$, but $\Delta^{\times}s\in S_+\tn
S_+$ by~\ref{lemma-delta*}. Elements of $k'$ may be moved between the
factors by the identifications 
$$
S'\bcp (\unt{k'}\bcp S'')=S'\bcp S''=(S'\bcp \unt{k'})\bcp S'', 
$$ 
and so the map descends to $\ct{S}\tn_{k'}\ct{S'}$.
Finally, it is a map of $k$-$k''$-bimodules by the
associativity of the composition product.

Since the map $\ct{S}\tn_{k'}\ct{S'}\to\ct{S\bcp_{k'}S'}$ is all we
need to make the cotangent space of a plethory into an algebra, we leave
the many details of the isomorphism to the reader. The key observation is that
$$s\lcp s'
=s\lcp(e+\ep^+(s'))\circ(e-\ep^+(s'))\circ s'
=s\circ(e+\ep^+(s'))\lcp (s'-\ep^+(s'))
$$
so that $S\bcp_{k'}S'$ is generated by elements of the form
$s\lcp s'$ with $s'\in S'_+$. 
This suggests the map of rings $f\:S\bcp S'\to k\oplus \ct{S}\tn\ct{S'}$
given by 
$f(s\lcp s')=\ep^+(s\lcp s')+(s\circ(e+\ep^+(s'))-\ep^+(s\lcp s'))\tn (s'-\ep^+(s'))$,
which descends to the inverse $\ct{S\bcp S'}\to\ct{S}\tn\ct{S'}$.
\end{proof}

\subsection{} \lb{sbsc-cotan-alg}
{\em $\ct{P}$ is an algebra over $k$.}  Let $P$ be a
$k$-plethory.  The composition $P\bcp P\to P$ and unit $\unt{k}\to P$
induce $\ct{P}\otimes_k\ct{P}\to \ct{P}$ and $k=\ct{\unt{k}}\to \ct{P}$
making $\ct{P}$ an algebra over $k$. Note that $e$ is the
unit for composition and thus the unit of this algebra.

\subsection{} \lb{sbsc-cotan-module}
{\em $I/I^2$ is a $\ct{P}$-module.}
Let $I$ be a $P$-ideal in a $P$-ring $R$.  Then by
\ref{prop-ideal-product}, $\ct{P}$ acts as a monoid on $I/I^2$.   
But~\ref{lemma-add-1st-order} implies this action is $\zz$-linear, and
we always have $(\alpha+\beta)\circ x = \alpha\circ x + \beta\circ x$;
so, this action is actually a $\ct{P}$-module structure on $I/I^2$.
The two $k$-module structures on $I/I^2$, one by way of $k\to\ct{P}$ and the
other $k\to R$, agree.

\section{Twisted bialgebras and their coactions}
\lb{sec-bialg}

First we recall some basic notions introduced by Sweedler~\cite{Sweedler:IHES},
as modified by Takeuchi~\cite[4.1]{Takeuchi}.

\subsection{}
If $A$ and $B$ are two algebras over $k$, then $A\tn_k B$, where the
$k$-module structure on each factor is given by multiplication on the
left, has two remaining $k$-actions: one by right multiplication on
$A$ and one by right multiplication on $B$.  Let $A\swe B$, the {\em
Sweedler product}, denote the subgroup where these two actions
coincide. It is an algebra over $k$ with multiplication
\[
\Big(\sum_i a_i\tn b_i\Big)\Big(\sum_ja'_j\tn b'_j\Big)=
\sum_{i,j} a_ia'_j\tn b_ib'_j.
\]
The Sweedler product is symmetric
in the sense that the symmetrizing map
\[
A\tn_k B \to B\tn_k A,\quad a\tn b\to b\tn a
\]
sends $A\swe B$ isomorphically to $B\swe A$.  Note that
$\swe$ is not naturally associative in the
generality above (but it is if, say, the algebras are 
$k$-flat~\cite[\S 2]{Sweedler:IHES}).

If $M$ and $N$ are left $A$-modules, then $M\tn_k N$ is a left
$A\swe A$-module by $(\sum_i a_i\tn b_i)(m\tn n)=\sum_i
a_im\tn b_in$.  

\subsection{}
\lb{sbsc-twisted}
We say $A$ is a {\em twisted $k$-bialgebra} if it is equipped with
a map $\Delta\:A \to A\swe A$ of algebras over $k$ and a map $\ep\:A \to k$
of $k$-modules
satisfying the following properties
\begin{enumerate}
\item the composite
$A\longlabelmap{\Delta} A\swe A \hookrightarrow A\tn A$
is coassociative with counit $\ep$, and
\item $\ep(1)=1$ and for all $a,b\in A$, we have 
$\ep(ab) = \ep(a\iota(\ep(b)))$, 
where $\iota$ denotes the structure map $k\to A$.
\end{enumerate}
Thus, the structure of a twisted $k$-bialgebra on $A$ is the same as
the structure of a $k$-bialgebroid on $A$ where the structure map
$k\tn_{\zz} k \to A$ factors through multiplication $k\tn_{\zz}
k\to k$.  (Several equivalent formulations of the notion of
bialgebroid are discussed in
Brzezinski--Militaru~\cite{Brzezinski-Militaru}.)  Assuming
flatness, it is also the same as what
Sweedler~\cite{Sweedler:IHES} called a $\times_k$-bialgebra structure.

The category of left $A$-modules then has a monoidal structure that is
compatible with $\tn_k$, and this is precisely the data needed to
make this so (\cite[5.1]{Schauenburg}\cite[3.1]{Brzezinski-Militaru}).
If $\Delta$ is cocommutative in the obvious sense, this monoidal
category is symmetric.

\subsection{}
\lb{sbsc-coaction}
Let $C$ be an algebra over $k$.  A coaction of $A$ on $C$
is a map $\alpha\:C\to A\swe C$ of algebras commuting with the maps
from $k$ such that the composite
\[
C\longlabelmap{\alpha} A\swe C \hookrightarrow A\tn C
\]
is a coaction of $A$, viewed as a $k$-coalgebra, on $C$.  (So, $C$ is
a left $A$-comodule algebra in the terminology of~\cite{Dascalescu:HA}).
Given a left $A$-module $M$, a left $C$-module
$N$, and a coaction of $A$ on $C$, the tensor product $M\tn_k N$
is naturally a left $C$-module by way of $\alpha$.
In this way, the category of left
$A$-modules acts on the category of left $C$-modules.

The map $\Delta\:A\to A\swe A$ is a coaction, the {\em regular}
coaction.

\subsection{} \lb{sbsc-semi-direct}
{\em Generalized semi-direct product $R\sd{A} C$.}
Suppose $A$ coacts on $C$ and also acts on a
$k$-ring $R$ in the sense that the multiplication map $R\tn
R\to R$ is a map of $A$-modules.  Then 
$R\tn_k C$
is an $R$-module and (by~\ref{sbsc-coaction}) a
$C$-module, and this induces a multiplication
\[
(R\tn C) \tn (R\tn C) = R\tn \big(C \tn (R\tn
C)\big) \longmap R\tn (R\tn C) \longmap R\tn C
\]
on $R\tn C$ with unit $1\tn 1$.  The map $k\to R\tn C$ is
simply $x\mapsto x\tn 1=1\tn x$.  

We denote this algebra by $R\sd{A} C$.
When $C$ is $A$ with the regular coaction and $A$ is untwisted (i.e., the
image of $k$ is in the center of $A$), this agrees with the semi-direct,
or ``smash'', product in the usual sense~\cite{Dascalescu:HA}.

It is immediate that the map $R\to R\sd{A} C$ given by $r\mapsto r\tn
1$ is a map of algebras over $k$, and the counit property
implies the map $C\to R\sd{A} C$, $c\mapsto 1\tn c$ is also such a
map.  Therefore an $R\sd{A} C$-module structure on a $k$-module $M$
is the same as actions of $R$ and $C$ on $M$ which are intertwined
as follows:
\[
c\big(r(c'm\big)) = \sum_i (c^{(1)}_i r)\big(c^{(2)}_i c' m\big),
\]
where $\Delta(c) = \sum_i c^{(1)}_i\tn c^{(2)}_i \in A\tn C$.

\section{The additive bialgebra}
\lb{sec-add-bialg}

The purpose of this section is to show that the set of additive
elements in a $k$-plethory is naturally a cocommutative twisted
$k$-bialgebra, at least under certain flatness hypotheses.

\subsection{}
\lb{sbsc-bialg-data}
Let $P$ be a $k$-plethory.  An element $f\in P$ is {\em additive} if
$\Delta^+(f)=f\tn 1 + 1\tn f$, which is equivalent to requiring that
$f(x+y)=f(x)+f(y)$ for all elements $x,y$ in all $P$-rings.  (In fact,
taking $x=e\tn 1$, $y=1\tn e$ in $P\tn P$ suffices).  Because we have
$\ep^+(f)=f(0)=0$, every additive element is in $P_+$. The set $A$, or
$A_P$, of
additive elements is clearly
closed under addition and composition, and composition by additive
elements distributes over addition; thus $A$ is a generally
non-commutative algebra with unit $1_A=e$.  Furthermore, the map
$\iota\:k\to A, c\mapsto ce$ is a map of algebras; so in this way, $A$ is
an algebra over $k$.

\begin{prop}
\lb{prop-interlinear}
The image of $A\swe A$ in $P\tn P$ is the set of $k$-interlinear
elements, where $f\in P\tn P$ is said to be $k$-interlinear if
$f(r,s)$ is additive in each argument $r,s\in R$ and we have
$f(cr,s)=f(r,cs)$ for all $c\in k$.
\end{prop}

Here we are using the notation of~\ref{sbsc-binary}.  Note that a
$k$-interlinear element $f$ is not required to be $k$-linear in each
argument.

\begin{proof}
First we show that the image of $A\tn P$ is the set of elements that
are additive on the left.  If $f$ is in the
image of $A\tn P$, it is immediate that $f$ is additive on the
left.  Now suppose $f$ is additive on the left.  Because $A$ is the
kernel of the $k$-module map
\bel{equ-beta}
P \longlabelmap{\varphi} P\tn P, \text{\ \ \ \ }
f \mapsto \Delta^+(f)-f\tn 1-1\tn f,
\end{equation}
the image of $A\tn P$ in $P\tn P$ is the kernel of the map 
$\varphi\tn 1\:P\tn P\to P\tn P\tn P$; so it is enough to show $f$ is in
the kernel of $\varphi\tn 1$.  Write $f=\sum_i\alpha_i\otimes\beta_i$.
Then we have
\begin{align*}
\sum_i\Delta^+(\alpha_i)\tn\beta_i &=
\sum_i\alpha_i(e\tn 1+1\tn e)\tn\beta_i \\
 &= \sum_i(\alpha_i\tn\beta_i)\circ(e\tn 1\tn 1+ 1\tn e\tn 1, 1\tn 1\tn e) \\
 &= f(e\tn 1\tn 1+ 1\tn e\tn 1, 1\tn 1\tn e) \\
 &= f(e\tn 1\tn 1, 1\tn 1\tn e) + f(1\tn e\tn 1, 1\tn 1\tn e) \\
 &= \sum_i \alpha_i(e)\tn 1\tn\beta_i(e) + 1\tn\alpha_i(e)\tn \beta_i(e) \\
 &= \sum_i (\alpha_i\tn 1+1\tn\alpha_i)\tn\beta_i. \\
\end{align*}
But $(\varphi\tn 1)(f)$ is the difference between the first and last
sums, and so $f$ is in the kernel.

Essentially the same argument shows 
the image of $A\tn A$ in $A\tn P$ is the set of
elements whose image in $P\tn P$ is both left additive and
right additive.

It is clear that any element in the image of $A\swe A$ is interlinear.
Now let $f$ be a $k$-interlinear element of $A\tn A$.
Then $f(ce\tn 1,1\tn e) = f(e\tn 1, 1\tn ce)$.
Writing $f=\sum_i \alpha_i\tn\beta_i$, we have
\[
\sum_i \big(\alpha_i\circ(ce)\big)\tn\beta_i 
   = \alpha_i\tn\big(\beta_i\circ(ce)\big),
\]
that is, $f$ transforms the same way under the two actions of $k$ on
$A$ by right multiplication.
\end{proof}

\begin{prop}
\lb{prop-bialg}
$\Delta^{\times}(A)$ is contained in the image of $A\swe A$
in $P\tn P$. If the maps $A^{\tn 2}\to P^{\tn 2}$ and
$A^{\tn 3}\to P^{\tn 3}$ are injective, the algebra $A$
is a cocommutative twisted $k$-bialgebra (\ref{sbsc-twisted}), where $\ep$ is
$\ep^{\times}$ and $\Delta$ is $\Delta^{\times}$, viewed as a map
$A\to A\swe A \subseteq A\tn A$.
\end{prop}

\begin{proof}
For any element $f\in A$ and any $P$-ring $R$, the map $R\times R\to
R$ given by $(r,s)\mapsto f(rs)$ is clearly $k$-interlinear.  Because
this map is just the application of $\Delta^{\times}(f)$, we see
$\Delta^{\times}(f)$ is $k$-interlinear and therefore
lies in the image of $A\swe A$, by~\ref{prop-interlinear}.

Now we show $\Delta$ is a map of algebras over $k$.  Take $a,b\in A$.
Because $a$ is additive and using~\ref{prop-P-equ}, we have
\[
\Delta^{\times}(a\circ b) = a\circ\Delta^{\times}(b) =
\sum_{i,j}(a^{[1]}_i\circ b^{[1]}_j) \tn (a^{[2]}_i\circ
b^{[2]}_j),
\]
but this last term is the product in $A\swe A$ of $\Delta^{\times}(a)$
and $\Delta^{\times}(b)$.  It is clear that $\Delta$ is a map over $k$.

The cocommutativity of $\Delta$ follows from that of $\Delta^{\times}$.

It remains to check properties (1)--(2) of~\ref{sbsc-twisted}.
Because we have $A\tn A\tn A \subseteq
P\tn P\tn P$, the coassociativity of $\Delta$ can be
tested in $P\tn P\tn P$, where it follows from the
associativity of the comultiplication $\Delta^{\times}$ on $P$.  
The map $\ep$ is a counit for $\Delta$ simply because $\ep^\times$ is
for $\Delta^\times$.

It is clear that $\ep(1)=1$.  By~\ref{prop-P-equ}, we also have
($\iota$ denoting the structure map $k\to A$)
\[
\ep^{\times}(a\circ \iota(\ep^{\times}(b)))
  = a\circ\ep^{\times}(e\ep^{\times}(b)) = a\circ((e\ep^{\times}(b))(1))
  = a\circ\ep^{\times}(b) = \ep^{\times}(a\circ b),
\]
for all $a,b\in A$.
\end{proof}

\subsection{} \lb{rmk-flatness}
{\em Remark.}
If $A$ and $P$ are flat over $k$, then the injectivity
hypotheses of the proposition hold.  In particular, they do if $k$ is
a Dedekind domain and $P$ is torsion-free.  They also hold if the
inclusion $A\to P$ is split, for example if $P=S(A)$.

In fact, we do not know any examples of plethories where the
assumptions of the previous proposition are not satisfied, but if they
exist, it seems clear that the correct replacement of $A$ would be the
collection of all multilinear elements in all tensor powers of $P$
assembled together in some sort of operadic coalgebra construction.

\section{The coaction of $A_P$ on $\ct{P}$}
\lb{sec-cotan2}

Because $A=A_P$ is contained in $P_+$, we have a map $A\to\ct{P}$, which
is clearly a map of algebras over $k$.

\begin{prop}
There is a unique map $\nu$ such that the diagram
\[
\xymatrix{
P_+\ar[rr]^-{\Delta^{\times}}\ar@{>>}[d]
  && P_+\tn P_+ \ar@{>>}[d] \\
\ct{P}\ar@{-->}[rr]^-{\nu}
  && P_+\tn\ct{P}  \\
}
\]
(using~\ref{lemma-delta*}) commutes; and the image of $\nu$ is
contained in the image of $A\swe\ct{P}$.

If the maps 
$A^{\tn i}\tn \ct{P} \to P_+^{\tn i}\tn\ct{P}$ are
injective for $i=1,2$, 
then $\nu$, viewed as a map $\ct{P}\to A\swe\ct{P}$, is
a coaction of the twisted $k$-bialgebra $A$
on $\ct{P}$, and the natural map $A\to\ct{P}$ is $A$-coequivariant, where
$A$ has the regular coaction.
\end{prop}

The injectivity hypotheses hold under the
flatness and splitting hypotheses of~\ref{rmk-flatness}.

\begin{proof}
The first statement is immediate because 
$\Delta^{\times}\:P\to P\tn P$ is a ring map. 

Let $\varphi$ be as in (\ref{equ-beta}).
To show the image of $\nu$ is contained in the image of
$A\tn\ct{P}$, it is enough to show the
composite map along the bottom row of the diagram
\[
\xymatrix{
P_+ \ar[r]^-{\Delta^\times}\ar@{>>}[d] 
  & P_+\tn P_+ \ar[r]^-{\varphi\tn 1}\ar@{>>}[d] 
  & P_+\tn P_+\tn P_+ \ar@{>>}[d] \\
\ct{P} \ar[r]^-{\nu} 
  & P_+\tn\ct{P}\ar[r]^-{\varphi\tn 1} 
  & P_+\tn P_+\tn\ct{P}
}
\]
is zero, and hence it is enough to show the composite of the maps along
the top and the right is zero.  The method is the same as that
of~\ref{prop-interlinear}. 

For any $f\in P_+$, write $\Delta^\times(f)=\sum_i f^{[1]}_i\tn
f^{[2]}_i$.  Then
\[
\sum_i \Delta^+(f^{[1]}_i)\tn f^{[2]}_i
  = \sum_i f^{[1]}_i(e\tn 1+1\tn e)\tn f^{[2]}_i(e) 
  = f((e\tn 1+1\tn e)\tn e).
\]
On the other hand, by \ref{lemma-add-1st-order} we can write
$\Delta^+(f) \equiv f\tn 1 + 1\tn f \mod J^2$, 
where $J=P\tn P_++P_+\tn P$.  Therefore we have
\begin{align*}
f(e\tn 1\tn e + 1\tn e\tn e) 
  &\equiv f(e\tn 1\tn e) + f(1\tn e\tn e)  
   \mod P\tn P\tn P_+^2\\
  &= \sum_i \left(f^{[1]}_i(e)\tn 1+1\tn f^{[1]}_i(e)
    \right)\tn f^{[2]}_i
\end{align*}
and hence
\begin{align*}
(\varphi\tn 1)(\Delta^{\times}(f)) 
  &= \sum_i \left(\Delta^+(f^{[1]}_i)-f^{[1]}_i(e)\tn 1+1\tn
                     f^{[1]}_i(e)\right)\tn f^{[2]}_i \\
  & \equiv 0 \mod  P\tn P\tn P_+^2,\\
\end{align*}
which was to be proved.

As in~\ref{prop-interlinear}, we show $\nu(f)$ is contained in the image
of $A\swe\ct{P}$ by applying $f$ to the equation $ce\tn e = e\tn ce$,
for any $c\in k$.

Now we show $\nu$ is a map of algebras over $k$.  Suppose
$f,g\in\ct{P}$, and write $\nu(f)=\sum_i f^{[1]}_i\tn f^{[2]}_i$
and $\nu(g)=\sum_j g^{[1]}_j\tn g^{[2]}_j$ with
$f^{[1]}_i,g^{[1]}_j\in A$ and $f^{[2]}_i,g^{[2]}_j\in\ct{P}$.  Then
\begin{align*}
\nu(f\circ g) &= f\left(\sum_j g^{[1]}_j\tn g^{[2]}_j\right)
\text{\quad by \ref{prop-P-equ}} \\ 
 &= \sum_j f\left(g^{[1]}_j\tn g^{[2]}_j\right) \\
 &= \sum_{i,j} \left(f^{[1]}_i\circ g^{[1]}_j\right) \tn \left(f^{[2]}_i\circ
      g^{[2]}_j\right) \\
 &= \left(\sum_i f^{[1]}_i\tn f^{[2]}_i\right)\circ
    \left(\sum_j g^{[1]}_j\tn g^{[2]}_j\right) \\
 &= \nu(f)\nu(g).
\end{align*}
And $\nu$ is a map over $k$ because
$\Delta^{\times}(ce)=c(e\tn e)$.

All that remains is to show that $\ep$ is a counit and that $\nu$ is
coassociative.  The first follows immediately from the counit property
of $\ep^{\times}$, and because of our assumptions, coassociativity can
be tested in $A\tn\ct{P}\tn\ct{P}$, where it follows from the fact
that $\Delta^{\times}$ is coassociative on $P$.
\end{proof}

\subsection{} \lb{sbsc-add-bialg}
{\em Example.}
If $B$ is a cocommutative $k$-bialgebra and $P=S(B)$, then
$\ct{P}=B$.  The image of inclusion $\ct{P}=B\hookrightarrow S(B)$ is
contained in $A$, and this is a section of the natural map
$A\to\ct{P}=B$.  The coaction of $A$ on $B$ is given by this
inclusion: 
\[
B\longlabelmap{\Delta} B\tn B \longmap A\tn B.
\]
If $k$ is a $\qq$-ring, the inclusion $B\hookrightarrow A$ is an
isomorphism, but if $k$ is an $\ff_p$-ring for some prime number
$p$, it will never be.  For we have $e^p\in A$, but the image of
$e^p$ in $\ct{P}$ is zero because $p\geq 2$.

\subsection{} 
{\em $I/I^2$ is an $R/I\sd{A}\ct{P}$-module.}
Let $I$ be a $P$-ideal in a $P$-ring $R$.  Then by
\ref{sbsc-cotan-alg}, $I/I^2$ is naturally a $\ct{P}$-module.
It follows from the associativity of the action of $P$ on $R$ that the
$\ct{P}$-action and $R/I$-action are intertwined as
in~\ref{sbsc-semi-direct}, and therefore these actions extend to an
action of $R/I\sd{A}\ct{P}$.

\subsection{} \lb{sbsc-action-diff}
{\em $\Omega^1_{R/k}$ is an $R\sd{A}\ct{P}$-module.}  
Let $R$ be a $P$-ring.  Because we have $\Omega^1_{R/k}=I/I^2$, where
$I$ is the kernel of the multiplication map $R\tn R\to R$, the
$R$-module $\Omega^1_{R/k}$ is naturally a $R\sd{A}\ct{P}$-module.

\section{Classical plethories revisited}
\lb{sec-construction}

Let $p$ be a prime number.
In this section we present a construction of $\Lambda_p$
(of~\ref{sbsc-lambda-p}), and hence an approach to the
$p$-typical Witt
vectors, which given the generalities developed earlier in this paper, is
completely effortless.  We also discuss the linearization of
$\Lambda_p$ and similar classical plethories.

\subsection{} 
\lb{sbsc-lambda-mod}
Consider the trivial $\ff_p$-plethory $\unt{\ff_p}$.  The bialgebra
$A$ of additive elements of $\unt{\ff_p}$ (see~\ref{prop-bialg}) is
the free bialgebra $\ff_p[F]$ on the monoid $\nn$ generated by the
Frobenius element $F=e^p$.  It therefore has a canonical lift $\zz[F]$
to a commutative bialgebra over $\zz$.  Let $\zz\ab{F}$ denote
$S(\zz[F])=\zz[F^{\circ\nn}]$, the free $\zz$-plethory on this
bialgebra.  The natural map $\zz\ab{F}\to\unt{\ff_p}$ is a surjection,
and so we can consider the amplification of $\zz\ab{F}$ along
$\unt{\ff_p}$.

\begin{prop}
There is a unique map of $\zz\ab{F}$-rings from $\Lambda_p$ to
the amplification of $\zz\ab{F}$ along $\unt{\ff_p}$, and this map is
an isomorphism of $\zz\ab{F}$-plethories.
\end{prop}

\begin{proof}
Let $P'$ denote the amplification.  Because $\Lambda_p$ is
$p$-torsion-free,~\ref{cor-amplification} implies we need only show
that a $\zz\ab{F}$-deformation of a $\unt{\ff_p}$-ring is the same as
a $p$-torsion-free $\Lambda_p$-ring.  But this is
just~\ref{lemma-cdd}, the strengthened form of Cartier's
Dieudonn\'e--Dwork lemma.
\end{proof}

\subsection{}
\lb{sbsc-frobenius-mod}
The same process gives ramified and twisted versions of the Witt
ring.  Let $\sO$ be a Dedekind domain, let $k$ be a residue field
of characteristic $p$, let $q$ be a power of $p$, and let $F$ be a
lift to $\sO$ of the endomorphism $x\mapsto x^q$ of $k$.  Then the
$\zz$-plethory $\zz\ab{F}$ acts on $\sO$, and we can form the plethory
$\sO\ab{F}\dfn\sO\tn \zz\ab{F}$, which maps to $\unt{k}$ by $F\mapsto
e^q$.  Let $M$ denote the rank-one $\sO$-module $\m^{-1}(F-e^q)$, and
let $B$ denote $\unt{\sO}\tn S_{\sO}(M)$.  One can easily check there
is a unique $\sO$-$\sO$-biring structure on $B$ such that the
inclusion $B\to K\ab{F}$ is a map of birings. (The structure maps are
similar to those in (\ref{equ-theta-coprod}).)  Let $P$ denote the free
pointed $\sO$-plethory on $B$.  Then an action of $P$ on an
$\m$-torsion-free $\sO$-ring $R$ is the
same~\ref{sbsc-free-on-pointed-biring}
as a map $B\bcp R\to R$ such that $e\lcp r \mapsto r$,
which is the same as an endomorphism $F\:R\to R$ extending the $F$ on
$\sO$ such that $F(x)\equiv x^q \mod \m$ for all $x\in R$.  Thus an
$\sO\ab{F}$-deformation of an $\unt{k}$-ring is the same as a $P$-action
on an $\m$-torsion-free $\sO$-ring.  Because $P$ is $\m$-torsion-free,
\ref{cor-amplification} gives a canonical isomorphism from $P$ to the
amplification of $\sO\ab{F}$ along $\unt{k}$.

When $\m$ is a principal ideal, surely much of this theory
agrees with Hazewinkel's formula-based approach~\cite[Ch.\
25]{Hazewinkel:Book} to objects of the same name.  Any precise results
along these lines would require some proficiency in his theory, which
proficiency we do not have.

It seems worth mentioning, however, that when $\m$ is not principal, 
it is unlikely $W_P(R)$ has
a description in terms of traditional-looking Witt components.  The
reason is simply that
the analogue of $W_2(R)$, the ring of length-two $\Lambda_p$-Witt vectors with
entries in $R$, is
\[
\Ring_{\sO}(B,R)=\Ring_{\sO}(\unt{\sO}\tn S_{\sO}(M),R)  
= R \times (\m\tn_{\sO} R),
\]
which is not naturally $R\times R$ (as sets).

\subsection{}
It is also possible to recover $\Lambda$ in this manner.  For a finite
set $S$ of prime numbers, construct a $\zz$-plethory $\Phi_S$ as follows: Let
$\Phi_{\{\}}$ denote the $\zz$-plethory $\zz\ab{\psi_p\,|\, p \text{
prime}}$, where the $\psi_p$ are ring-like (\ref{eg-bialg}) and
commute with each other.  For $S'=S\cup \{p\}$, where $p$ is a
prime not in $S$, let $\Phi_{S'}$ denote the
amplification of $\Phi_S$ along $(\ff_p\tn\Phi_S)/(\psi_p^{\circ
n}-e^{p^n}\,|\, n\geq 0)$.  Using induction, one can construct a natural
map $\Phi_S\to\Lambda$ and prove that $\Phi_S$ is
torsion-free and that torsion-free $\Phi_{S}$-rings are
the same as torsion-free $\Phi_{\{\}}$-rings such that
$\psi_p(x)\equiv x^p\mod p$ for all $p\in S$.  It is also possible to show
that $\Phi_S$ is canonically independent of the order of the amplifications.

Using Wilkerson's result~\cite{Wilkerson}
that a torsion-free $\lambda$-ring is the same a
ring equipped with commuting Adams operations $\psi_p$ such that
$\psi_p(x)\equiv x^p\mod p$ for all primes $p$,
it follows that the maps
$\Phi_S\to\Lambda$ induce an isomorphism from the colimit of the
$\Phi_S$ to $\Lambda$.

One could certainly construct variants for rings of
integers in general number fields, as in the single-prime case above.

\subsection{} 
{\em Linearization of $\Lambda_p$.} The additive bialgebra of
$\Lambda_p$ is $\zz[F]$ with comultiplication $F\mapsto F\tn F$.
(Because $\Lambda_p$ is torsion-free, additivity can be checked in
$\qq\otimes\Lambda_p=\qq\ab{F}$, to which~\ref{sbsc-add-bialg} can be
applied.)  It follows---either from the traditional, explicit
description (\ref{sbsc-lambda-p}) of $\Lambda_p$ or
from~\ref{sbsc-free-on-pointed-biring}---that the cotangent space
$\ct{\Lambda_p}$ is freely generated by the image $\theta$ of
$\theta_1$, the coaction is given by $\theta\mapsto F\tn
\theta$, and the map $\zz[F]\to\zz[\theta]$ is $F\mapsto
p\theta$.  Note that
\[
\theta^n = p^{-n}F^n \equiv \theta_n \mod(\Lambda_p)_+,
\]
that is, the two familiar generating sets $\{\theta_n\}$ and
$\{\theta_1^{\circ n}\}$ of $\Lambda_p$ agree in
$\ct{\Lambda_p}$.  Also note that the map $F\mapsto \theta$ is an
isomorphism from $A$ to $\ct{\Lambda_p}$ of algebras with an
$A$-coaction, but the canonical map is not this map, or even an
isomorphism at all. The general case of~\ref{sbsc-frobenius-mod} is very
similar, but there is no canonical element $\theta$, only 
$\m^{-1}F$.

\subsection{}
{\em Linearization of $\Lambda$.}
The situation for $\Lambda$ is essentially the same.  Its additive bialgebra is
$\zz[\psi_p\,|\, p\text{ prime}]$ with $\Delta\:\psi_p\mapsto\psi_p\tn\psi_p$.
The cotangent space is $\zz[\lambda_p\,|\, p\text{ prime}]$, and the
coaction of $\zz[\psi_p\,|\, p\text{ prime}]$ on $\zz[\lambda_p\,|\, p\text{
prime}]$ is given by $\lambda_p\mapsto\psi_p\tn\lambda_p$.  The map
$\zz[\psi_p\,|\, p\text{ prime}]\to\zz[\lambda_p\,|\, p\text{ prime}]$ is given
by $\psi_p\mapsto (-1)^pp\lambda_p$.  These can be checked using
Newton's formulas~\cite[I (2.11)$'$]{MacDonald:SF}.

\subsection{}
The binomial plethory is $\Lambda/(\psi_n-e\mid n\geq 1)$;
its additive bialgebra is the trivial one, $\zz$, and its
cotangent algebra is $\qq$.

\subsection{} 
{\em Bloch's Frobenius.} There is an endomorphism of the de~Rham--Witt
complex~\cite{Illusie:dRW}, which is usually called Frobenius, but which on $i$-forms is
$p^{-i}F$, where $F$ is the actual Frobenius map.  In fact, this
endomorphism lifts to the de~Rham complex of $W(R)$:
By~\ref{sbsc-action-diff}, the element $\theta\in\ct{P}$ acts on
$\Omega^1_{W(R)}$, but we have $\theta=p^{-1}F\in\ct{P}$, and so
$\theta$ reduces to Bloch's Frobenius map in degree $1$.  In degree
$i>0$, Bloch's Frobenius is $\theta^{\tn i}$
as in
\[
\theta^{\tn i}(\eta_1\wedge\cdots\wedge\eta_i) =
\theta(\eta_1)\wedge\cdots\wedge\theta(\eta_i). 
\]

We remarked above that there is an isomorphism of $A_{\Lambda_p}$ and
$\ct{\Lambda_p}$ of algebras with an $A_{\Lambda_p}$-coaction
identifying $F$ and $\theta$ but that this is not the
canonical map.  This is perhaps a pleasant explanation of the meaning of the
well-known fact that Bloch's is {\em a} Frobenius
operator even though it is not {\em the} Frobenius operator.

For the variant of $\Lambda_p$ over a general integer ring $\sO$ at a prime
$\m$, the compatibility between any
generalization of Bloch's Frobenius map and the true one would involve
some choice of uniformizer, and so it would be a mistake to try to
find such a generalization.  
Instead it is the $\sO$-line $\m^{-i}F^{\tn i}=(\m^{-1}F)^{\tn i}$
that acts.

\subsection{}
{\em Remark.}  The perfect closure $(\unt{\ff_p})^{p^{-\infty}}$ of
the ring $\unt{\ff_p}$ has a unique $\ff_p$-plethory structure
compatible with that of $\unt{\ff_p}$.  Let $\zz\ab{F^{\circ\pm 1}}$
denote the free $\zz$-plethory on the group bialgebra $\zz[F,F^{-1}]$
of $\zz$.  Then the map of plethories $\zz\ab{F^{\circ\pm
1}}\to(\unt{\ff_p})^{p^{-\infty}}$ is a surjection.  One can show the
amplification $P$ of this map is the plethory push-out, or amalgamated
product, of $\Lambda_p$ and $\zz\ab{F^{\circ\pm 1}}$ over $\zz\ab{F}$.
Its Witt functor is particularly interesting and useful: if $V$ is an
$\ff_p$-ring, $W_P(V)$ is $A_{\mathrm{inf}}(V/\zz_p)$, the universal
$p$-adic formal pro-infinitesimal $\zz_p$-thickening of $V$, in the
sense of Fontaine~\cite[1.2]{Fontaine:p-adic-periods}.


\bibliography{references}

\begin{thebibliography}{10}

\bibitem{AT:Lambda}
M.~F. Atiyah and D.~O. Tall.
\newblock Group representations, {$\lambda $}-rings and the {$J$}-homomorphism.
\newblock {\em Topology}, 8:253--297, 1969.

\bibitem{Bergman}
George~M. Bergman and Adam~O. Hausknecht.
\newblock {\em Co-groups and co-rings in categories of associative rings},
  volume~45 of {\em Mathematical Surveys and Monographs}.
\newblock American Mathematical Society, Providence, RI, 1996.

\bibitem{Brzezinski-Militaru}
Tomasz Brzezi{\'n}ski and Gigel Militaru.
\newblock Bialgebroids, {$\times\sb A$}-bialgebras and duality.
\newblock {\em J. Algebra}, 251(1):279--294, 2002.

\bibitem{Cartier:Witt-gen}
Pierre Cartier.
\newblock Groupes formels associ\'es aux anneaux de {W}itt g\'en\'eralis\'es.
\newblock {\em C. R. Acad. Sci. Paris S\'er. A-B}, 265:A49--A52, 1967.

\bibitem{Dascalescu:HA}
Sorin D{\u{a}}sc{\u{a}}lescu, Constantin N{\u{a}}st{\u{a}}sescu, and
  {\c{S}}erban Raianu.
\newblock {\em Hopf algebras}, volume 235 of {\em Monographs and Textbooks in
  Pure and Applied Mathematics}.
\newblock Marcel Dekker Inc., New York, 2001.
\newblock An introduction.

\bibitem{Deligne:Ex-Ser}
Pierre Deligne.
\newblock La s\'erie exceptionnelle de groupes de {L}ie.
\newblock {\em C. R. Acad. Sci. Paris S\'er. I Math.}, 322(4):321--326, 1996.

\bibitem{Fontaine:p-adic-periods}
Jean-Marc Fontaine.
\newblock Le corps des p\'eriodes {$p$}-adiques.
\newblock {\em Ast\'erisque}, (223):59--111, 1994.
\newblock With an appendix by Pierre Colmez, P\'eriodes $p$-adiques
  (Bures-sur-Yvette, 1988).

\bibitem{Hazewinkel:Book}
Michiel Hazewinkel.
\newblock {\em Formal groups and applications}, volume~78 of {\em Pure and
  Applied Mathematics}.
\newblock Academic Press Inc. [Harcourt Brace Jovanovich Publishers], New York,
  1978.

\bibitem{Illusie:dRW}
Luc Illusie.
\newblock Complexe de de\thinspace {R}ham-{W}itt et cohomologie cristalline.
\newblock {\em Ann. Sci. \'Ecole Norm. Sup. (4)}, 12(4):501--661, 1979.

\bibitem{Joyal:Witt}
Andr{\'e} Joyal.
\newblock {$\delta$}-anneaux et vecteurs de {W}itt.
\newblock {\em C. R. Math. Rep. Acad. Sci. Canada}, 7(3):177--182, 1985.

\bibitem{Knutson:Lambda}
Donald Knutson.
\newblock {\em {$\lambda $}-rings and the representation theory of the
  symmetric group}.
\newblock Springer-Verlag, Berlin, 1973.
\newblock Lecture Notes in Mathematics, Vol. 308.

\bibitem{Lang:Algebra-I}
Serge Lang.
\newblock {\em Algebra}.
\newblock Addison-Wesley Publishing Co., Inc., Reading, Mass., 1965.

\bibitem{Lazard:CFG}
Michel Lazard.
\newblock {\em Commutative formal groups}.
\newblock Springer-Verlag, Berlin, 1975.
\newblock Lecture Notes in Mathematics, Vol. 443.

\bibitem{MacLane:CWM}
Saunders Mac~Lane.
\newblock {\em Categories for the working mathematician}.
\newblock Springer-Verlag, New York, second edition, 1998.

\bibitem{MacDonald:SF}
I.~G. Macdonald.
\newblock {\em Symmetric functions and {H}all polynomials}.
\newblock Oxford Mathematical Monographs. The Clarendon Press Oxford University
  Press, New York, second edition, 1995.
\newblock With contributions by A. Zelevinsky, Oxford Science Publications.

\bibitem{Schauenburg}
Peter Schauenburg.
\newblock Bialgebras over noncommutative rings and a structure theorem for
  {H}opf bimodules.
\newblock {\em Appl. Categ. Structures}, 6(2):193--222, 1998.

\bibitem{Sweedler:IHES}
Moss~E. Sweedler.
\newblock Groups of simple algebras.
\newblock {\em Inst. Hautes \'Etudes Sci. Publ. Math.}, (44):79--189, 1974.

\bibitem{Takeuchi}
Mitsuhiro Takeuchi.
\newblock Groups of algebras over {$A\otimes \overline A$}.
\newblock {\em J. Math. Soc. Japan}, 29(3):459--492, 1977.

\bibitem{Tall-Wraith}
D.~O. Tall and G.~C. Wraith.
\newblock Representable functors and operations on rings.
\newblock {\em Proc. London Math. Soc. (3)}, 20:619--643, 1970.

\bibitem{Wilkerson}
Clarence Wilkerson.
\newblock Lambda-rings, binomial domains, and vector bundles over {${\bf
  C}P(\infty )$}.
\newblock {\em Comm. Algebra}, 10(3):311--328, 1982.

\bibitem{Witt:Vectors}
Ernst Witt.
\newblock Zyklische {K}\"orper und {A}lgebren der {C}harakteristik $p$ vom
  {G}rad $p^n$. {S}truktur diskret bewerter perfekter {K}\"orper mit
  vollkommenem {R}estklassen-k\"orper der charakteristik $p$.
\newblock {\em J. Reine Angew. Math.}, (176), 1937.

\bibitem{Witt:Works}
Ernst Witt.
\newblock {\em Collected papers. {G}esammelte {A}bhandlungen}.
\newblock Springer-Verlag, Berlin, 1998.
\newblock With an essay by G\"unter Harder on Witt vectors, Edited and with a
  preface in English and German by Ina Kersten.

\bibitem{Wraith}
G.~C. Wraith.
\newblock Algebras over theories.
\newblock {\em Colloq. Math.}, 23:181--190, 325, 1971.

\end{thebibliography}
\bibliographystyle{plain}
\end{document}